\input amstex
\voffset=-8mm
\documentstyle{amsppt}
\NoBlackBoxes
\magnification=1200
\parindent 20 pt

\define\1{^{-1}}

\define\CP{\Bbb C \Bbb P}
\define\N{\Bbb N}

\define \C{\Bbb C}
\define \Z{\Bbb Z}

\define \pr{\operatorname{pr}}
\define \Cl{\operatorname{Cl}}

\define \Int{\operatorname{Int}}

\define \tm—

\define\F{\Sigma}
\define\Si{\frak S}
\define\R{F}
\define\r{\varepsilon}

{\catcode`\@=11
\gdef\n@te#1#2{\leavevmode\vadjust{%
 {\setbox\z@\hbox to\z@{\strut#1}%
  \setbox\z@\hbox{\raise\dp\strutbox\box\z@}\ht\z@=\z@\dp\z@=\z@%
  #2\box\z@}}}
\gdef\leftnote#1{\n@te{\hss#1\quad}{}}
\gdef\rightnote#1{\n@te{\quad\kern-\leftskip#1\hss}{\moveright\hsize}}
\gdef\?{\FN@\qumark}
\gdef\qumark{\ifx\next"\DN@"##1"{\leftnote{\rm##1}}\else
 \DN@{\leftnote{\rm??}}\fi{\rm??}\next@}}

\def\refACC   {1}
\def\refBa   {2}
\def\refBT   {3}
\def\refBK   {4}
\def\refBr  {5}
\def\refCh   {6}
\def\refChi   {7}
\def\refG   {8}
\def\refGr   {9}
\def\refKK   {10}
\def\refKT   {11}
\def\refLM  {12}
\def\refL  {13}
\def\refMP   {14}
\def\refMW   {15}
\def\refMo   {16}
\def\refMoi   {17}
\def\refMT   {18}
\def\refO   {19}
\def\refR   {20}
\def\refSi   {21}
\def\refSik   {22}
\def\refV   {23}

\topmatter
\title On Braid Monodromy Factorizations \endtitle
\author V. Kharlamov and Vik. S. Kulikov \endauthor

\centerline {\it
To Igor Rostislavovich
Shafarevich
on his 80th birthday} \vskip1cm

\abstract We introduce and develop a language of semigroups over
the braid groups for a study of braid monodromy factorizations
(bmf's) of plane algebraic curves and other related objects. As an
application we give a new proof of Orevkov's theorem on
realization of a bmf over a disc by algebraic curves and show that
the complexity of such a realization can not be bounded in terms
of the types of the factors of the bmf. Besides, we prove that the
type of a bmf is distinguishing Hurwitz curves with singularities
of inseparable types up to $H$-isotopy and $J$-holomorphic
cuspidal curves in $\C P^2$ up to symplectic isotopy.

\endabstract
\thanks \flushpar The first author is a member of Research Training
Networks EDGE and RAAG, supported by the European Human Potential
Program. The second author is partially supported by
INTAS-00-00259, 00-00269 and RFBR  02-01-00786.
\endthanks
\endtopmatter

\document

\baselineskip 20pt
\subheading{Introduction}

In this paper, we deal with algebraic curves and other related
objects in the projective plane, such as $J$-holomorphic and
Hurwitz curves (the definition of Hurwitz curves is given in
Section 3), which imitate the behavior of plane algebraic curves
with respect to pencils of lines. A common feature of all these
geometric objects is that for each of them there can be defined so
called braid mondromy factorizations (bmf's, in short), which are
known to be a powerful tool for study the topology of embedding of
curves in $\C P^2$. Since foundating works by O.~Chisini
\cite{\refCh},\cite{\refChi}, these braid monodromy factorizations
are considered as genuine factorizations in the braid groups and
studied up to some moves, called in our days Hurwitz moves. We
propose to study them by means of suitable semigroups over the
braid groups, so that Hurwitz equivalent factorizations become
represented by elements of these semigroups. In some cases (such
as that of topological Hurwitz curves, see Section 3) it is useful
to go up to a second level, i.e., to consider factorization
semigroups over semigroups of the first level.

As it seems to us, the language of semigroups simplifies
constructing and study of the objects defined by bmf's. As an
example we give an almost pure algebraic proof of a recent result
due to S.~Orevkov \cite{\refO}, which states that any bmf over a
disc can be realized by an algebraic curve (and which generalizes
Rudolf's theorem \cite{\refR} on algebraic realization of
quasi-positive braids). In our construction the curve seats in a
ruled surface and has only the simplest ramifications outside the
disc. For arbitrary, non necessary round, discs the construction
is explicit so that the degree of the curve and the ruling can be
bounded. On the other hand, we show, by means of Moishezon
examples \cite{\refMoi}, that for round discs there does not exist
any bound in terms of the types of the factors of the bmf.

A braid monodromy factorization of a projective curve in a ruled
surface is a factorization of $\Delta^{2N}$, where $\Delta$ is the
so called Garside element and $N$ is the degree of the ruling. Its
factors correspond to the critical values of the restriction of
the projection to the curve. As is known, contrary to over disc
factorizations not any projective bmf can be realized by an
algebraic curve (see \cite{\refMoi}). On the other hand, any
projective bmf can be realized in the class of Hurwitz curves. We
prove that the type of a braid monodromy factorization (by the
type we mean the orbit under the natural conjugacy action of the
braid group, see Section 1) is distinguishing Hurwitz curves with
singularities of types $w^k=z^n$ up to $H$-isotopy (i.e., isotopy
in the class of Hurwitz curves) at least in the case when the
critical values in distinct critical points are distinct. And we
show that any such Hurwitz curve is $H$-isotopic to an
almost-algebraic one, i.e., a one which can be given by an
algebraic equation over a disc containing all the critical values
of the projection.

In Section 4 we give few remarks on isotopies of $J$-holomorphic
curves in $\C P^2$. In particular, using known results on
$J$-holomorphic curves and the results of Section 3, we show how
various symplectic isotopy problems are reduced to a pure
algebraic, in a sense, study of braid factorization types. In this
section we show that two cuspidal $J$-holomorphic curves are
symplectically isotopic if and only if they have the same  braid
monodromy factorization type.
We also give a bmf-characterization of algebraicity for nodal
symplectic surfaces in $\C P^2$: a nodal symplectic surface is
symplectically isotopic to an algebraic curve  if and only if its
bmf is a partial re-degeneration of a factorization whose factors
are conjugates of the squares of standard generators of the braid
group.

The paper is organized as follows. Section 1 is devoted to
factorization semigroups: we give basic constructions, investigate
their functorial properties, and apply them to the braid groups.
There we introduce the notion of stable equivalence in
semigroups over a braid group and prove that two elements with
conjugated factors are stably equivalent
if an only if they factorize the same
element of the braid group. This result is then used in Section 2
for the proof of the generalized Rudolf theorem mentioned above.
Section 3 is devoted to $H$-istopies. There, besides Hurwitz and
almost-algebraic curves, which have algebraic singularities, we
consider what we call topological Hurwitz curves allowing them to
have arbitrary cone singularities.
Then, we introduce a class of cone singularities of inseparable type and
prove that two topological Hurwitz curves of the same degree
having singularities of inseparable types
any two of which lie in different fibers of the projection
are $H$-isotopic if and only if these curves
have the same braid monodromy factorization type.
Section 4 deals with the symplectic case.

{\bf Acknowledgements.} We are grateful to E.~Artal-Bartolo,
S.~Nemirovski and V.~Schevchishin for usefull discussions and proposals.
This research was started during the stay of the second author
in Strasbourg university and finished within the frame
of RiP programs in Mathematisches Forschungsinstitut Oberwolfach.

\subheading{\S 1.\ Factorization calculus}

{\bf 1.1. Factorization semigroups.}
A collection $(S,B,\alpha,\lambda) $,
where $S$ is a semigroup, $B$ is a group, and $\alpha :S\to B$,
$\lambda :B\to \text{Aut}(S)$ are homomorphisms, is
called {\it a semigroup $S$ over a group $B$} if for all
$s_1,s_2\in S$
$$s_1\cdot s_2=\lambda
(\alpha (s_1))(s_2)\cdot s_1= s_2\cdot\rho (\alpha (s_2))(s_1),$$
where $\rho (g)=\lambda(g^{-1})$.
If we are given two semigroups $(S_1,B_1,\alpha_1,\lambda_1)$
and $(S_2,B_2,$ $\alpha_2,\lambda_2)$
over, respectively, groups $B_1$ and $B_2$,
we call a pair $(h_1,h_2)$ of homomorphisms $h_1:S_1\to S_2$ and
$h_2: B_1\to B_2$
{\it a homomorphism of semigroups over
groups} if
\roster \item"{($i$)}" $h_2\circ \alpha_{S_1}=\alpha_{S_2}\circ h_1$,
\item"{($ii$)}" $\lambda_{B_2}(h_2(g))(h_1(s))=h_1(\lambda_{B_1}(g))(s)$
for all $s\in S_1$ and all $g\in B_1$.
\endroster

The {\it factorization semigroups} defined below constitute the
principal, for our purpose, examples of semigroups over groups.

Let $\{ g_i \}_{i\in I}$ be a set of elements of a group $B$. For
each $i\in I$ denote by $O_{g_i}\subset B$ the set of all the
elements in $B$ conjugated to $g_i$ (the orbit of $g_i$ under the
action of $B$ by inner automorphisms). Call their union $X=
\cup_{i\in I} O_{g_i}\subset B$ the {\it full set of conjugates}
of $\{ g_i\}_{i\in I}$ and the pair $(B,X)$ an {\it equipped group}.

For any full set of conjugates $X$ there are two natural maps
$r=r_X: X\times X \to X$ and  $l=l_X: X\times X \to X$  defined by
$ r(a,b)=b^{-1}ab$ and  $ l(a,b)=aba^{-1}$ respectively. For each
pair of letters $a,b\in X$ denote by $R_{a,b;r}$ and
$R_{a,b;l}$ the relations defined in the following way:
$$
R_{a,b;r}\quad
\text{stands for} \quad a\cdot b=b\cdot\,r(a,b)\quad\text{if $b\ne\bold 1$
and} \,\, a\cdot\bold 1=a\quad\text{otherwise};
$$
$$
R_{a,b;l} \quad
\text{stands for} \quad a\cdot b=l(a,b)\cdot\,a \quad\text{if $a\ne\bold 1$
and} \,\, \bold 1\cdot b=b\quad\text{otherwise}.
$$
Now, put
$$
\Cal{R}=\{ R_{a,b;r},
 R_{a,b;l} \,\vert \, (a,b)\in X\times X,\, a\ne b\, \, \text{if}\,\,
a\ne \bold 1 \, \,
\text{or}\,\, b\ne \bold 1 \}
$$
and introduce
the semigroup
$$S(B,X)= \langle\, x\in X\,\,
:\,\, R\in \Cal{R}\, \rangle 
$$
by means of this relation set
$\Cal R$.
Introduce also a homomorphism $\alpha_X :S(B,X)\to B$
given by $\alpha_X (x)=x$ for each $x\in X$.

Next, we define two actions $\lambda$ and $\rho$ of the group $B$
on the set $X$:
$$ x\in X \mapsto \rho (g) (x)=g^{-1}xg\in X $$
and
$$ x\in X\mapsto \lambda (g) (x)=gxg^{-1}\in X.$$

As is easy to see, the above relation set $\Cal{R}$ is preserved
by the both actions and, therefore, $\rho$ and $\lambda$ define an
anti-homomorphism $\rho : B\to \text{Aut} (S(B,X))$ ({\it right
action}) and a homomorphism $\lambda : B\to \text{Aut} (S(B,X))$
({\it left  or conjugation action}). The action $ \lambda (g)$ on
$S(B,X)$ is called {\it simultaneous conjugation } by $g$. Put
$\lambda _{S} =\lambda \circ\alpha_X $ and $\rho _{S} =\rho
\circ\alpha_X$.

\proclaim{Claim 1.1} For any $g\in B$ and any $x_i,x_j\in X$ we
have \roster \item"{($i$)}" $\lambda (g)=\rho (g^{-1})$;
\item"{($ii$)}" $\alpha_X (r(x_i,x_j))=x_j^{-1}x_ix_j$;
\item"{($iii$)}" $\alpha_X (l(x_i,x_j)=x_ix_jx_i^{-1}$;
\item"{($iv$)}" $\rho _{S} (x_i)(x_j)=r(x_j,x_i)$; \item"{($v$)}"
$\lambda _{S} (x_i)(x_j)=l(x_i,x_j)$; \item"{($vi$)}" $\rho
(\alpha_X (x_i)^{-1})(x_j)=l(x_i,x_j)$; \item"{($vii$)}" $\lambda
(\alpha_X (x_i)^{-1})(x_j)=r(x_j,x_i)$.
\endroster
\endproclaim

\demo{Proof} Straightforward. \qed\enddemo

It follows from Claim 1.1 that $(S(B,X),B, \alpha_X,
\lambda_S)$ is a semigroup over $B$. We call such semigroups the
{\it factorization semigroups over} $B$. When $B$ is fixed, we
abbreviate $S(B,X)$ to $S_X$. By $x_{1}\cdot\dots\cdot x_{n}$ we
denote the element in $S_X$ defined by a word $x_{1}\dots x_{n}$.

Notice that $S: (B,X)\mapsto (S(B,X),B, \alpha_X, \lambda)$
is a functor from the category of
equipped groups to the category of the semigroups over groups.
In particular, if $X\subset Y$ are
two full sets of conjugates in $B$, then
the identity map $id: B\to B$ defines
an embedding $id_{X,Y}:S(B,X)\to S(B,Y)$.
So that, for each group $B$, the semigroup $S_B=S(B,B)$
is an {\it universal factorization semigroup over} $B$,
which means that each semi-group $S_X$ over $B$
is canonically embedded in $S_B$ by $id_{X,B}$.

Since $\alpha_X=\alpha_B \circ id_{X,B}$,
we make no difference between $\alpha_X$ and $\alpha_B$
and denote the both simply by $\alpha$.

Denote by $B_{X}$ the subgroup of $B$ generated by
the image of $\alpha :S(B,X)\to B$, and for each
$s\in S_X$ denote by $B_s$
the subgroup of $B$ generated by the images
$\alpha (x_1),\dots ,\alpha (x_n)$ of the
elements $x_1, \dots , x_n$
of a factorization $s=x_1\cdot\dots \cdot x_n.$

\proclaim{Claim 1.2} The subgroup $B_s$ of $B$
does not depend on the presentation of
$s$ as a word in letters $x_i$ of $X$.
\endproclaim

\demo{Proof} It follows from ($ii$) and ($iii$) of Claim 1.1.
\qed\enddemo

\proclaim{Proposition 1.1}
For any $X$ and any $s\in S_X$ as above,
\roster
\item"{($i$)}" $\text{ker}\, \lambda $ coincides with
the centralizer $C_{X}$ of $B_{X}$ in $B$;
\item"{($ii$)}" if $\alpha (s)$ belongs to the center $C(B_s)$
of $B_s$, then $\lambda (g)$ leaves fixed $s\in S_X$
whatever is $g\in B_s$.
\endroster
\endproclaim

\demo{Proof} ($i$) is evident.

($ii$) The group $B_s$
is generated in $B$ by $\alpha(x_1),\dots ,\alpha(x_n)$,
where $s=x_1\cdot\dots \cdot x_n$ with $ x_i\in X$. Therefore,
to prove ($ii$) it is sufficient to show that
$\lambda_S(x_i)(s)=s$ for each $i=1,\dots ,n$
as soon as $\alpha(s)\in C(B_s)$.
Using the relations $x_j\cdot x_i=x_i\cdot r(x_j,x_i)$, we can move
$x_i$ to the left and obtain a presentation of $s$ in the form
$$s=x_i\cdot\widetilde x_1\cdot\dots\cdot \widetilde x_{n-1}=
x_i\cdot\widetilde s.
$$
If $\alpha(s)\in C(B_s)$, then
$$
\lambda_S(s)(x_i)=x_i.
$$
Finally,
$$
\align
s & =l(\widetilde x_1,x_i)\cdot \dots\cdot
l(\widetilde x_{n-1},x_i)\cdot x_i=\lambda_S(x_i)
(\widetilde x_1\cdot \dots\cdot\widetilde x_{n-1})\cdot x_i= \\
& =\lambda_S(x_i)(\widetilde s)\cdot x_i=
\lambda_S(\lambda_S(x_i)(\widetilde s))(x_i)\cdot\lambda_S(x_i)(\widetilde s)=
\lambda_S(x_i\cdot\widetilde s )(x_i)\cdot\lambda_S(x_i)(\widetilde s)=\\
& =x_i\cdot\lambda_S(x_i)(\widetilde s)=
\lambda_S(x_i)(x_i)\cdot\lambda_S(x_i)(\widetilde s)=\lambda_S(x_i)(s).
\endalign
$$
\qed\enddemo

Consider two full sets of conjugates $X_1,X_2$ in $B$ and the
semigroups $S_{X_1}$ and $S_{X_2}$ associated with them. A map
$\psi : X_2\to S_{X_1}$ can be extended to a homomorphism
$\psi : S_{X_2}\to S_{X_1}$ if and only if for all
$x_i,x_j\in X_2$ the equalities
$$\psi(x_i)\cdot\psi(x_j)=\psi(x_j)\cdot\psi(r(x_i,x_j))$$
and
$$\psi(x_i)\cdot\psi(x_j)=\psi(l(x_i,x_j))\cdot\psi(x_i)$$
hold in $S_{X_1}$.

We say that the homomorphism $\psi$ is defined over $B$
if $\alpha_{X_2}(x)=\alpha_{X_1}(\psi(x))$ for all $x\in X_2$.
\newline
{\bf Example 1.1.} Let $X_1$ be the set of the conjugates of an
element $x_1\in B$ and $X_2$ the set of the conjugates of $x_1^2$.
Assume that the map $\phi :X_1\to X_2$ given by $\phi(x)=x^2$ for
$x\in X_1$ is a bijection. Then the map $\psi :X_2\to S_{X_1}$
given by $\psi(x)=\phi ^{-1}(x)\cdot \phi^{-1}(x)$ defines a
homomorphism  $\psi :S_{X_2}\to S_{X_1}$ over $B$.
\newline
{\bf Example 1.2.} Example 1.1 can be generalized as follows. Pick
a $n$-set $\{x_1,\dots ,x_n\}$ of elements in a full conjugate set
$X_1\subset B$ and a $k$-set of products $s_j(x_1,\dots ,
x_n)=x_{i_1(j)} \cdot\dots\cdot x_{i_{m(j)}(j)}\in S_{X_1}$,
$j=1,\dots ,k$. Consider $X_2=O_{\bar s_1}\cup\dots \cup O_{\bar
s_k}$ with $O_{\bar s_j}$ being the full set of conjugates of
$\bar s_j=\alpha(s_j)\in B$. Assume that $\bar s_i$ and $\bar s_j$
are not conjugated in $B$ for $i\neq j$. Then, the map $X_2\to
S_{X_1}$ given by $g\bar s_jg^{-1} \mapsto \lambda(g)(s_j)\in
S_{X_1}$ can be extended uniquely to a homomorphism $r:S_{X_2}\to
S_{X_1}$ defined over $B$. Such a homomorphism $r$ is called {\it
re-degeneration} of the set $\{ s_j\}$.

In Section 4 we use a kind of generalization of this notion (which
no more takes a form of a homomorphism). It looks as follows. In
notation of Example 1.2, put $Z=X_1\cup X_2$ and consider an
element $z=z_1\cdot z_2\in S_Z$ where $z_1\in S_{X_2}$, $z_2\in
S_Z$. The element
$$\overline z=r(z_1)\cdot z_2\in S_Z$$
is called a {\it partial re-degeneration} of $z$.

The construction of the semigroups $S(B,X)$ can be iterated.
Namely, one can consider the conjugation action
of $B$ on $S(B,X)$, pick any set $Y$ which is a union
of orbits of this action and
introduce the semigroup $S(S(B,X),Y)$
as a semigroup generated by the letters $s\in Y$
and being subject to the relations
$$ s_i\overlineú s_j=s_j\overline{\cdot }\rho _S(s_j)(s_i)$$
and
$$s_i\overlineú s_j=\lambda_S(s_i)(s_j)\overlineú s_i$$
for all $s_i, s_j\in Y$.

One can introduce, in addition, the homomorphisms
$\beta_S: S(S(B,X),Y)\to S(B,X)$ sending
$s=s_1\overlineú\dots\overlineú s_n\in S(S(B,X),Y)$ to
$s_1\cdot\dots\cdot s_n\in S(B,X)$
and the conjugation actions $\lambda : B \to
\text{Aut}(S(S(B,X),Y))$ ($\lambda (g)$ is
acting as {\it simultaneous conjugation} by $g$),
as well as associated with them
the homomorphisms $\beta=\alpha\circ \beta_S$
and the actions
 $\lambda_{S}=\lambda \circ\alpha_S :
S(B,X) \to \text{Aut}(S(S(B,X),Y))$,
 $\lambda_{S,S}=\lambda \circ \beta :
S(S(B,X), Y) \to \text{Aut}(S(S(B,X),Y))$
of, respectively, $S(B,X)$ and $S(S(B,X),Y)$ on
$S(S(B,X),Y)$.
The right actions $\rho$,
$\rho_S$, and $\rho_{S,S}$ of, respectively,
$B$, $S(B,X)$, and $S(S(B,X), Y)$ on $S(S(B,X),Y)$
are defined in a similar way.

On the other hand, if $X$ is a subset of $Y$, then there is the
natural embedding of $S_X=S(B,X)$ into $S(S_X,Y)$. Moreover, there
is the natural embedding of $S(S_X,Y)$ into the universal
semigroup $S(S_B)=S(S_B,S_B)$ over $S_B$. Thus, any semigroup
$S_X$ over $B$ can be considered as a subsemigroup of $S(S_B)$ and
we may, without any confusion, denote the operation
$\overlineú$ in $S(S_B)$ by $\cdot$. Note that $S(S_X,X)$ is
naturally isomorphic to $S_X$.

{\bf 1.2. Hurwitz equivalence.} Let, as above,
$Y$ be a union of orbits of the
conjugation action $\lambda$ of $B$ on $S(B,X)$.
An ordered set
$$\{ y_1,\dots ,y_n \,\,\mid y_i\in Y \}, n\in\Z$$
is called a { \it factorization} of
$g=\beta (y_1)\dots\beta (y_n)\in B$ in $Y$.
Denote by $F_{X,Y}\subset \bigcup_n Y^n$ the set of all
possible factorizations of the elements of $B$ in $Y$
over all $n\in \N$.
There is a natural map $\varphi : F_{X,Y} \to S(S(B,X),Y)$, given by
$$\varphi(\{ y_1,\dots ,y_n \})=y_1\cdot\dots \cdot y_n .$$
The transformations which
replace in $\{ y_1,\dots ,y_n \}$ some two neighboring factors
$(y_i, y_{i+1})$ by $(y_{i+1},\rho_S(y_j)(y_{i+1})$
or
$(\lambda_S(y_i)(y_{i+1}),y_i)$ and preserve the other factors
are called {\it Hurwitz moves}. Two factorizations are
{\it Hurwitz equivalent} if one can be obtained
from the other by a finite sequence of Hurwitz moves.

\proclaim{Claim 1.3} Two factorizations $y=\{ y_1,\dots ,y_n \}$ and
$z=\{ z_1,\dots ,z_n \}$ are Hurwitz equivalent if and only if
$\varphi (y)=\varphi (z)$.
\endproclaim

\demo{Proof} Evident. \qed\enddemo
{\bf Remark 1.1.} Below, according with Claim 1.3,
we identify classes of
Hurwitz equivalent factorizations in $Y$ with their images in $S(S_X,Y)$.
And when $Y=X$, we identify $S(S_X,Y)$ with $S_X$.

{\bf 1.3. Semigroups over the braid group
and stable equivalence.} In this subsection, $B=B_m$ is the braid group
with $m$ strings. We fix a set $\{ a_1,\dots ,a_{m-1} \} $ of
so called {\it standard generators}, i.e.,
generators being subject to the relations
$$\align
a_ia_{i+1}a_i & =a_{i+1}a_i a_{i+1} \qquad \qquad 1\leq i\leq n-1 , \tag 1.1 \\
a_ia_{k} & =a_{k}a_i  \qquad \qquad \qquad \, \, \mid i-k\mid \, \geq 2.
\tag 1.2
\endalign
$$
We denote by $B_m^+$ the semi-group
defined by the same generating letters and relations.

\proclaim{Garside's Theorem} \cite{\refG} The natural homomorphism
$i:B_m^+\to B_m$ is an embedding.
\endproclaim

Following this theorem,
we identify $B_m^+$ with its image $i(B_m^+)$ in
$B_m$ and call the images $i(g)$ of the elements $g\in B_m^+$
{\it positive} elements of the group $B_m$.

Denote by $A_k=A_k(m), \,k\ge 0,$ the full
set of conjugates of $a_1^{k+1}$ in $B_m$ (recall that all the generators
$a_1,\dots ,a_{m-1}$ are conjugated to each other).
Consider the semigroup $S_{A_0}$ as a subsemigroup of
the universal semi-group $S_{B_m}$ over $B_m$. A positive word
$g=a_{i_1}\dots a_{i_n}$ in the alphabet $\{a_1,\dots, a_n\}$
defines an element
$\overline g(a_1,\dots ,a_{m-1})=a_{i_1}\cdot\dots\cdot a_{i_n}\in S_{A_0}$.
On the other hand, $g$ defines an element $\widetilde g=
a_{i_1}\dots a_{i_n}$ of $B_m^+$.

\proclaim{Lemma 1.1} \cite{\refBT} A  map $\nu : B^+_m \to
S_{A_0}$ given by $\nu(\widetilde g)=\overline g$ is a
well-defined injective homomorphism of semi-groups.
\endproclaim

\demo{Proof} To show that $\nu$ is a well-defined
homomorphism, it is sufficient to check that
the relations (1.1) and (1.2) hold in $S_{A_0}$. We have
$$
\align
a_i\cdot a_{i+1}\cdot a_i &
=a_{i+1}\cdot (a_{i+1}^{-1}a_ia_{i+1})\cdot a_i
 = a_{i+1}\cdot a_i\cdot (a_i^{-1}a_{i+1}^{-1}a_ia_{i+1}a_i)= \\
& =a_{i+1}\cdot a_i\cdot (a_i^{-1}a_{i+1}^{-1}a_{i+1}a_{i}a_{i+1})
= a_{i+1}\cdot a_i\cdot a_{i+1}
\endalign
$$
for $1\le i\le n-1$ and
$$
a_i\cdot a_k= a_k\cdot(a_k^{-1}a_ia_k)=a_k\cdot a_i
$$
for $\vert i-k\vert\ge 2$.
The homomorphism $\nu$ is injective, since $\alpha_B \circ \nu$ is the identity
map by Garside's Theorem.
\qed\enddemo

Let $\Delta =\Delta _m $ be
the so-called Garside element:
$$
\Delta=
(a_1\dots a_{m-1})\dots(a_1a_2a_3)(a_1a_2)a_1.
$$
As is well-known,
$$ \Delta^2=(a_1\dots a_{m-1})^{m}$$
is the generator of the center of $B_m$.
Denote by $\delta ^2=\delta _m^2$ the element
in $S_{A_0}\subset S_{B_m}$ equal to
$$ \delta^2=(a_1\cdot ...\cdot a_{m-1})^{m}.$$

\proclaim{Lemma 1.2} The element $\delta ^2$ is
fixed under the conjugation action of $B_m$ on $S_{B_m}$,
i.e., $\rho (g)(\delta ^2)=\delta ^2$ for any $g\in B_m$.
\endproclaim

\demo{Proof} It follows from
$\alpha(\delta^2)=\Delta^2$ and Proposition 1.1 ($ii$)
applied to $s=\delta^2$ (for which $(B_m)_s=B_m$).
\qed\enddemo

In our study of topological Hurwutz surfaces (see section 3.2)
we use an extention $\widetilde S_{B_m}$ of 
$S_{B_m}$ which is defined as follows. To each element $I$
belonging to the set $\Cal I$ of all the subsets of $\{1,\dots ,
m\}$, let
us associate a letter $\bold 1_{I}$. 
Consider a semigroup $\widetilde S_{B_m}$ generated by the pairs
$(g,\bold 1_I)$, $g\in B_m$ and  $I\in \Cal I$, and being subject
to the relations
$$(g_1,\bold 1_{I_1})\cdot (g_2, \bold 1_{I_2})=
(g_1g_2g_1^{-1},\bold 1_{\sigma (g_1)(I_2)})\cdot (g_1,\bold
1_{I_1}), $$
$$ (g_1,\bold 1_{I_1})\cdot (g_2, \bold 1_{I_2})=
(g_2, \bold 1_{I_2})\cdot(g_2^{-1}g_1g_2,\bold 1_{\sigma
(g_2^{-1})(I_1)})$$ for all $g_1, g_2\in B_m$ and all $I_1,I_2\in \Cal I$;
here $\sigma$ is an action on $\Cal I$ induced by the
natural homomorphism from $B_m$ to the symmetric group $\Sigma _m$
acting on $\{1,\dots ,m\}$.
To extend the actions $\lambda$ and $\rho$ of $B_m$ on $S_{B_m}$
to the action on $\widetilde S_{B_m}$ we put
$$\lambda (b)((g,\bold 1_I))=(\lambda(b)(g),\bold 1_{\sigma(b)(I)})
$$
and $\rho (g)=\lambda (g^{-1})$ for $b,g\in B_m$ and $I\in \Cal I$.
Also we extend the homomorphism $\alpha $ by
$$\alpha ((g,\bold 1_I))=g\in B_m$$
for all $I$. Note that the map $g\mapsto (g,\bold 1_{\emptyset})$
is extended to an embedding of $S_{B_m}$ in $\widetilde S_{B_m}$
over $B_m$.

Denote by $B_{k,i}$, $k+i\leq m$, a subgroup of the braid group
$B_m$ generated by a part $a_{i+1},\dots, a_{i+k-1}$ of a fixed
set of standard generators $a_1,\dots , a_{m-1}$ of $B_m$. We say
that an element $b\in B_m$ has {\it the interlacing number}
$l(b)=k$ if $k$ is the smallest number such that $b$ is conjugated
in $B_m$ to an element in $B_{k,0}$. For a pair $(B_{n,i}, b)$
with $b\in B_{n,i}$ and $l(b)=k$, an element $\widetilde g=(\bar
b, \bold 1_{\{i+k+1,\dots , i+n\}})\in \widetilde S_{ B_m}$ is
called a {\it standard tbmf-form} of $b\in B_{n,i}$ if $\bar b\in
B_{k,i}\subset B_{n,i}$ is conjugated to $b$ in $B_m$ (if $i+k\geq
n$ then ${\{i+k+1,\dots , i+n\}}=\emptyset$).

Now to each finite sequence of integers $k_1,\dots,k_t$ such that
$k_1+\dots +k_t\leq m$, $k_1\geq 2, \dots, k_t\geq 2$, let
associate a sequence of subgroups $B_{k_i,k_1+\dots +k_{i-1}}$,
$1\leq i\leq t,$ of $B_m$. Introduce also the subsets
$T_{k_1,\dots ,k_{t}}$ of $\widetilde S_{B_m}$ consisting of the
products $s=\widetilde g_1\cdot\dots\cdot \widetilde g_t$, where
$\widetilde g_i\in \widetilde S_{B_m}$ are the standard tbmf-forms
of elements in $B_{k_i,k_1+\dots +k_{i-1}}$ for each $1\leq i\leq
t$. Then, define $T=T_m$ to be the union $\bigcup \lambda (g)
T_{k_1,\dots ,k_{t}}$ over all $g\in B_m$ and all sequences of
integers $k_1,\dots,k_t$ such that $k_1+\dots +k_t\leq m$. Note
that for any permutation $\sigma\in \Sigma _t$ and any $\widetilde
g_1\cdot\dots\cdot \widetilde g_t\in T_{k_i,k_1+\dots +k_{i-1}}$,
one has
$$\widetilde g_1\cdot\dots\cdot \widetilde g_t=
\widetilde g_{\sigma(1)}\cdot\dots\cdot \widetilde g_{\sigma(t)}.$$

The elements of $T$ are called {\it tbm factorizations} and the
semigroup $\Cal T=\Cal T_m= S(\widetilde S_{B_m},T)$ is called the
{\it tbm factorization semigroup}. Two tbm factorizations are said
of the same {\it factorization type} if they belong to the same
orbit under the conjugation action of $B_m$.

The group $B_m$ as the set can be represented as the disjoint union over
$k$, $1\leq k\leq m$, of the  orbits of $T_{k,\emptyset}=
\{ (b,\bold 1_{\emptyset})\mid l(b)=k\}\subset T_k$ under the conjugation
action of $B_m$.
This presentation defines an imbedding $i:S_{B_m}\to \Cal T$
of semigroups over $B_m$. Thus, when it can not lead to a
confusion, we identify $S_{B_m}$ with its image
$i(S_{B_m})\subset\Cal T$.

We say that an element $s_1\in \Cal T$
(in particular, $s_1\in S_{B_m}$)
is {\it stably equal} to an element $s_2\in \Cal T$
($s_2\in S_{B_m}$, respectively)
if there is an integer $n\ge 1$ such that
$$ s_1\cdot (\delta^2)^n=s_2\cdot (\delta^2)^n$$
in $\Cal T$.

\proclaim{Theorem 1.1} Let $O_{x_1}, \dots , O_{x_n}$
be the orbits of elements
$x_1, \dots , x_n\in T$
under the conjugation action of $B_m$ on $T$.
Then for any $y_i\in O_{x_i}$, $1\leq i\leq n$, and
for any permutation $\sigma \in \Sigma_n$, the elements
$s_1=x_1 \cdot \dots \cdot x_n$ and $s_2=y_{\sigma(1)} \cdot \dots \cdot
y_{\sigma(n)}$ are stably equal in $\Cal T$ if and only if
$\beta (s_1)=\beta (s_2)$.
\endproclaim

{\bf Remark 1.2.} The above theorem remains true if $T$ is replaced
by any sets $X\subset S_{B_m}$, containing $A_0$
and $\Cal T$ by $S(S_{B_m},X)$.

\demo{Proof of Theorem 1.1} It is evident that if $s_1$ and
$s_2$ are stably equal in $\Cal T$ then
$\beta (s_1)=\beta (s_2)$.

Let $\beta (s_1)=\beta (s_2)$.
Since any permutation is a product of transpositions,
we can assume that $\sigma = \text{id}$.
Indeed, for each $g_1, g_2\in \Cal T$ we have the relation
$g_1\cdot g_2=g_2\cdot \rho (g_2)(g_1)$ as a relation in
$\Cal T$ and
in which $g_1$ and $\rho (g_2)(g_1)$ belong to the same orbit.
So, applying these relations we
get $s_2=\tilde y_{1} \cdot \dots \cdot
\tilde y_{n}$ with the factors $\tilde y_{i}\in O_{x_i}$.

\proclaim{Lemma 1.3} (\cite{\refG}) For any $g\in B_m$ there are
positive elements $r_1,r_2\in B_m$ and integers $k, p\in \Z$,
$p\geq 1$, such that \roster
\item"{($i$)}" $g= \Delta^{2k}r_1$;
 \item"{($ii$)}" $gr_2 = \Delta^{2p}$.
 \endroster
\endproclaim

\demo{Proof} It follows from Theorem 5 in \cite{\refG}.
\qed\enddemo

By Lemma 1.3 ($i$), since $x_i$ and $\tilde y_{i}$ belong to the same
orbit $O_{x_i}$ and
$\Delta^{2}$ belongs to the center of $B_m$,
we may assume that there are positive elements $g_i$ such that
$\tilde y_{i}=\rho (g_i^{-1})(x_i)$.
Applying Lemma 1.3 ($ii$) to each $g_i$,
we can find positive elements $r_i$ and positive integers $p_i$ such that
$g_ir_i = \Delta^{2p_i}$. By Garside's Theorem, Claim 1.1 and Lemma 1.1,
$\bar g_i\cdot \bar r_i = (\delta^{2})^{p_i}$
in $S_{A_0}\subset S_{B_m}\subset \Cal T$.
Put $p=p_1+\dots +p_n$. By Proposition 1.1 ($ii$),
it follows that for each
$x\in S_{B_m}$, we have $x\cdot \delta^2=\delta^2 \cdot x$ in $S_{B_m}$.
In addition,
$\rho (g_i^{-1})(x_i)\cdot \bar g_i=\bar g_i \cdot x_i$.
Therefore,
$$\align
s_2\cdot (\delta^2)^p & =
\rho (g_1^{-1})(x_1)\cdot\dots\cdot \rho (g_n^{-1})(x_n)
\cdot (\delta^2)^p = \\
& =\rho(g_1^{-1})(x_1)\cdot (\delta^2)^{p_1} \cdot\dots\cdot
\rho (g_n^{-1})(x_n)\cdot (\delta^2)^{p_n} = \\
& =\rho (g_1^{-1})(x_1)\cdot \bar g_1\cdot \bar r_1
\cdot\dots\cdot \rho (g_n^{-1})(x_n)
\cdot \bar g_n\cdot \bar r_n = \\
& = \bar g_1 \cdot x_1\cdot \bar r_1
\cdot\dots\cdot \bar g_n\cdot x_n\cdot \bar r_n.
\endalign
 $$

Consider in the beginning the case when all $x_i$ are standard generators.
Then, applying the Garside theorem to
$\alpha(s_1)(\Delta^2)^p= \alpha(s_2)(\Delta^2)^p$ in $B_m$
we get
$$
s_1\cdot (\delta^2)^p=\bar g_1 \cdot x_1\cdot \bar r_1
\cdot\dots\cdot \bar g_n\cdot x_n\cdot \bar r_n
$$
in $\nu (B^+_m)$, which gives
$s_1\cdot (\delta^2)^p=s_2\cdot (\delta^2)^p$ in $\nu (B^+_m)$.

In general case, all $\bar g_i$ and $\bar r_i$ belong to $\nu (B^+_m)$.
Applying the relations $a_i\cdot x_j =x_j\cdot \rho_{S,S}(x_j)(a_i)$,
we can move to the left all $x_i$ and obtain that
$s_2\cdot (\delta^2)^p=s_1\cdot s_3$, where
$s_3=\prod (t_i^{-1}z_it_i)$ and each $z_i$ is a letter of the alphabet
$\{ a_1, \dots , a_{m-1} \}$. Thus there is a positive integer $q$
such that $s_3\cdot (\delta^2)^q=(\prod z_i)\cdot (\delta^2)^q$,
and therefore
$$s_2\cdot (\delta^2)^{p+q}=s_2\cdot (\delta^2)^p\cdot (\delta^2)^q
=s_1\cdot s_3\cdot (\delta^2)^q=s_1\cdot (\prod z_i)\cdot (\delta^2)^q=
s_1\cdot (\delta^2)^{p+q},$$
since $(\delta^2)^{p+q}$ and $(\prod z_i)\cdot (\delta^2)^q$ belong to
$\nu (B^+_m)$
and $\alpha((\prod z_i)\cdot (\delta^2)^q)=(\Delta^2)^{p+q}$.
\qed\enddemo

Let us address two problems which seem to be open.

{\bf Garside problem.} Is $\alpha : S_{A_0(m)}\to B_m$ an embedding for
any $m$? In particular, does the equation $\alpha(s)=\Delta_m^2$
have only one solution, $s=\delta_m^2$?

{\bf Word problem.} Does the word problem for $\Cal T_m$ (
respectively, for $S_{B_m}$,
$S_{A_{\le 2}}$ with $A_{\le 2}=\bigcup_{k\le 2}A_k(m)$)
have the positive solution?

{\bf 1.4. "Pure nodal" semigroup.}
In this subsection we work with the semigroup $S_{A_1(m)}$.
Let fix a set of standard generators $\{ a_1, \dots,
a_{m-1}\}$ of $B_m$ and consider $B_{m-1}$ as a subgroup of $B_m$ generated
by $\{ a_1, \dots, a_{m-2}\}$. Put
$$\tilde\delta_m^2=\prod_{l=m}^{2} \prod^{l-1}_{k=1} z_{k,l}^2 \in S_{A_1},$$
where $z_{k,l}=(a_{l-1}\dots a_{k+1})a_k(a_{l-1}\dots
a_{k+1})^{-1}$ for $k<l$ (the notation $\prod_a^b$ states for the
left to the right product from $a$ to $b$). As is known, see for
example \cite{\refMT}, $\alpha(\tilde\delta_m^2)=\Delta_m^2$ and
$$\tilde\delta_m^2=\prod^{m-1}_{k=1} z_{k,m}^2\cdot\tilde\delta_{m-1}^2.$$

\proclaim{Lemma 1.4}
\roster
\item"{($i$)}"
$a_kz_{j,k}^2a_k^{-1}=z_{j,k+1}^2$;
\item"{($ii$)}" $a_kz_{k,m}^2a_k^{-1}=z_{k+1,m}^2$;
\item"{($iii$)}" $a_kz_{j,r}^2a_k^{-1}=z_{j,r}^2$
if either $r<k$, or $j<k$, or $j>k+1$;
\item"{($iv$)}" $a_kz_{k+1,m}^2a_k^{-1}=
z_{k+1,m}^{-2}z_{k,m}^2z_{k+1,m}^{2}$.
\endroster
\endproclaim

\demo{Proof} It follows from relations (1.1) and (1.2) and the
definition of the elements $z^2_{k,j}$.
\qed\enddemo

\proclaim{Proposition 1.2}
The element $\tilde\delta_m^2 \in S_{A_1}$ is a fixed element under the
conjugation action of $B_m$.
\endproclaim

\demo{Proof} It is sufficient to show that
$\lambda (a_i) (\tilde\delta_m^2)=\tilde\delta_m^2$ for
$i=1,\dots, m-1$. These equalities will be proved by induction on $m$.

By induction hypothesis and Lemma 1.4, we have that for $i<m-1$
$$
\align
& \lambda (a_i) (\tilde\delta_m^2) =
\lambda (a_i)(\prod^{m-1}_{k=1} z_{k,m}^2\cdot\tilde\delta_{m-1}^2)=  \\
& =\prod_{k=1}^{i-1}\lambda (a_i)(z^2_{k,m})\cdot
\lambda (a_i)(z^2_{i,m})\cdot
\lambda (a_i)(z^2_{i+1,m})\cdot \prod_{k=i+2}^{m-1}
\lambda (a_i)(z^2_{k,m})
\cdot \lambda (a_i)(\tilde\delta_{m-1}^2)= \\
& =\prod_{k=1}^{i-1}z^2_{k,m}\cdot z^2_{i+1,m}\cdot
(z^{-2}_{i+1,m}z^{2}_{i,m}z^{2}_{i+1,m})\cdot \prod_{k=i+2}^{m-1}
z^2_{k,m} \cdot \tilde\delta_{m-1}^2  =\tilde\delta_m^2.
\endalign
$$

Applying again Lemma 1.4, we have that for $i=m-1$
$$
\align
& \lambda (a_{m-1}) (\tilde\delta_m^2) =
\lambda (a_{m-1})(\prod^{m-1}_{k=1} z_{k,m}^2\cdot
\prod^{m-2}_{k=1} z_{k,m-1}^2
\cdot\tilde\delta_{m-2}^2)=  \\
& =\prod_{k=1}^{m-2}\lambda (a_{m-1})(z^2_{k,m})\cdot
\lambda (a_{m-1})(z^2_{m-1,m})\cdot
\prod^{m-2}_{k=1}\lambda (a_{m-1})(z^2_{k,m-1})\cdot
\lambda (a_{m-1})(\tilde\delta_{m-2}^2)= \\
& =\prod_{k=1}^{m-2} (a_{m-1}z^2_{k,m}a_{m-1}^{-1})\cdot
z^2_{m-1,m}\cdot
\prod^{m-2}_{k=1}(a_{m-1}z^2_{k,m-1}a_{m-1}^{-1})\cdot
\tilde\delta_{m-2}^2= \\
& = z^2_{m-1,m}\cdot \prod^{m-2}_{k=1}(a_{m-1}^{-1}z^2_{k,m}a_{m-1})\cdot
\prod^{m-2}_{k=1}z^2_{k,m}\cdot
\tilde\delta_{m-2}^2= \\
& = z^2_{m-1,m}\cdot \prod^{m-2}_{k=1}z^2_{k,m-1}\cdot
\prod^{m-2}_{k=1}z^2_{k,m}\cdot
\tilde\delta_{m-2}^2.
\endalign
$$

To complete the proof of the Proposition, it is sufficient to show that
$$\prod^{m-2}_{k=1} z_{k,m}^2\cdot z^2_{m-1,m}\cdot
\prod^{m-2}_{k=1} z_{k,m-1}^2=
z^2_{m-1,m}\cdot \prod^{m-2}_{k=1}z^2_{k,m-1}\cdot
\prod^{m-2}_{k=1}z^2_{k,m}.
\tag 1.3 $$
We have
$$t_m=\alpha (z^2_{m-1,m}\cdot \prod^{m-2}_{k=1}z^2_{k,m-1})=
a_{m-1}^2a_{m-2}\dots a_2a_1^2a_2\dots a_{m-2}
 $$
and to prove equality (1.3)
it is sufficient to show that
$$t_mz_{k,m}^2t_m^{-1}=z_{k,m}^2$$
for $k=1,\dots m-2$.
By induction on $m$, applying relations (1.1) and (1.2), we have
for $k\leq m-3$
$$
\align
& t_mz_{k,m}^2t_m^{-1}=(a_{m-1}^2a_{m-2}^{-1}t_{m-1}a_{m-2})z_{k,m}^2
(a_{m-1}^2a_{m-2}^{-1}t_{m-1}a_{m-2})^{-1}= \\
& =(a_{m-1}^2a_{m-2}^{-1}t_{m-1}a_{m-2}a_{m-1}a_{m-2})z_{k,m-2}^2
(a_{m-1}^2a_{m-2}^{-1}t_{m-1}a_{m-2}a_{m-1}a_{m-2})^{-1}= \\
& =(a_{m-1}^2a_{m-2}^{-1}t_{m-1}a_{m-1}a_{m-2}a_{m-1})z_{k,m-2}^2
(a_{m-1}^2a_{m-2}^{-1}t_{m-1}a_{m-1}a_{m-2}a_{m-1})^{-1}= \\
& =(a_{m-1}^2a_{m-2}^{-1}t_{m-1}a_{m-1}a_{m-2})z_{k,m-2}^2
(a_{m-1}^2a_{m-2}^{-1}t_{m-1}a_{m-1}a_{m-2})^{-1}= \\
& =(a_{m-1}^2a_{m-2}^{-1}t_{m-1}a_{m-1})z_{k,m-1}^2
(a_{m-1}^2a_{m-2}^{-1}t_{m-1}a_{m-1})^{-1}= \\
& =(a_{m-1}^2a_{m-2}a_{m-1}a_{m-2}^{-2}t_{m-1})z_{k,m-1}^2
(a_{m-1}^2a_{m-2}a_{m-1}a_{m-2}^{-2}t_{m-1})^{-1}= \\
& =(a_{m-1}a_{m-1}a_{m-2}a_{m-1}a_{m-2}^{-2})z_{k,m-1}^2
(a_{m-1}a_{m-1}a_{m-2}a_{m-1}a_{m-2}^{-2})^{-1}= \\
& =(a_{m-1}a_{m-2}a_{m-1}a_{m-2}^{-1})z_{k,m-1}^2
(a_{m-1}a_{m-1}a_{m-2}a_{m-1}a_{m-2}^{-1})^{-1}= \\
& =(a_{m-2}a_{m-1})z_{k,m-1}^2
(a_{m-2}a_{m-1})^{-1}=(a_{m-2}a_{m-1}a_{m-2})z_{k,m-2}^2
(a_{m-2}a_{m-1}a_{m-2})^{-1}= \\
& =(a_{m-1}a_{m-2}a_{m-1})z_{k,m-2}^2
(a_{m-1}a_{m-2}a_{m-1})^{-1}=a_{m-1}z_{k,m-1}^2a_{m-1}^{-1}=z_{k,m}^2.
\endalign
$$

Using the same calculations as above, one can show that
$t_mz^2_{m-2,m}t_m^{-1}=z^2_{m-2,m}$.
\qed\enddemo

The following theorem is a consequence of Proposition 1.2 and
Corollary 4.1, Theorem 3.1, Remark 4.1 which will be proven in
sections 3 and 4.

\proclaim{Theorem 1.2} Let $A_1\subset B_m$ be the full set of conjugates
of $a_1^2$. Then,
$\tilde\delta_m^2$ is the only element $s\in S_{A_1}$
such that $\alpha (s)=\Delta_m^2$. \qed
\endproclaim

\subheading{\S 2.\ Existence of polynomials with given braid monodromy
factorizations over a disc}

{\bf 2.1. Local braid monodromy over a point (unique germ).} Here
and further we denote by $(z,w)$ the standard coordinates in
$\C^2$. The proofs of the (well known) statements used here and
the further references can be found, for example, in
\cite{\refBK}.

If a germ $(C,0)\subset (\C^2,0)$ of a reduced complex analytic
curve does not contain the germ $z=0$
then it is given by an equation
$$P(z,w)=0, \tag 2.1$$
where
$$
P(z,w)=w^k+\sum_{i=1}^{k}q_{i}(z)w^{k-i},
$$
$q_i(z)$ are convergent power series (i.e., $q_i(z)\in\C\{z\}$),
$q_i(0)=0$ and the polynomial $w^k+\sum q_{i}(z)w^{k-i}\in \C\{
z\} [w]$ has no multiple factors. Therefore, one can choose a
small polydisc $D=D_1\times D_2\subset\C^2$,
$D_1=D_{1}(\varepsilon_1)=\{ z\in \C \, \mid \,\, \mid z\mid <
\varepsilon_1 \, \}$ and  $D_2=D_{2}(\varepsilon_2)= \{ w\in \C \,
\mid \,\, \mid w\mid < \varepsilon_2 \, \},$ such that: $C$ is an
analytic set at each point of the closure $\Cl D$ of $D$, the
projection on the $z$-factor $\pr=\pr_1: C\cap D\to D_{1}$ is a
proper finite map of degree $k$, and $(z,w)=0$ is the unique
critical point of $\pr_{\mid C\cap\Cl D}$. Reciprocally, if $D$ is
a polydisc and $C$ is $W$-prepared in $D$, i.e., if $C$ is a
reduced complex analytic curve with the latter properties with
respect to $\pr$, then it is defined in $\Cl D$ by an equation of
the same type. A $W$-prepared germ $(C,o)$ is {\it algebraic} (in
coordinates $(z,w)$) if and only if $q_i\in\C[z]$ for each
$i=1,\dots, k$.

To define the braid monodromy, let pick a point $u\in\partial D_1$
and put $D_{2,u}=\pr^{-1}(u)$,
$K=K(u)=\{w_1,\dots,w_k\}=D_{2,u}\cap C$. The loop $\partial D_1$
oriented counter-clockwise and starting at $u$ lifts to $\partial
D\cap C$ as a motion $\pr_2(\{w_1(t),\dots,w_k(t)\})$ of $k$
distinct points in $D_2$ starting and ending at $K$. This motion
defines a braid $b_{(C,o)}\in B_k=B_k[D_2,K]$ which is called the
{\it braid monodromy of $(C,o)$ with respect to} $\pr$. Note that
$l(b_{(C,o)})=k$.

The link of $(C,0)$ is an iterated positive torus link. This link
is determined by the Puiseaux pairs of irreducible components of
$(C,0)$ and the mutual intersection numbers of the components.
Therefore, as soon as we choose in $B_k$ as standard generators
the half-twists $a_1,\dots, a_{k-1}$, the braid $b_{(C,o)}$
becomes an element of  $B^+_k \subset B_k $ (a positive braid)
and is called the {\it standard form of the braid monodromy of
$(C,o)$ with respect to} $\pr$. The standard generators and the
standard form $b_{(C,0)}$ are {\bf defined uniquely up to
conjugation}.

The topological type of the triple $(D,C,\pr)$ is determined by
the standard form of its braid monodromy, and vice versa. Besides,
for each triple $(D,C,\pr)$ there is a constant $M=M_{(C,o)}\in
\N$ such that the topological type of $(D,C,\pr)$ coincides with
the topological type of a singularity given by
$\overline{P}(z,w)=w^k+\sum \bar q_{i}(z)w^{k-i}=0$, where $\bar
q_i(z)=q_i(z)+z^Mr_i(z)$ and $r_i(z)$ are arbitrary analytic
functions. In particular, one can make 
$C$ algebraic without changing the topological type of
$(D,C,\pr)$.

The topological type of the triple
$(D,C,\pr)$, and thus the standard form of
its braid monodromy, is determined by
a {\it resolution of singularities relative to $\pr$}.
By the latter we mean a sequence of
blow ups $\sigma_1 :
U_1\to D$, $\dots , \sigma_n:U_n\to U_{n-1}$
with centers at points
such that $\sigma^{-1}(C\cup \R)$,
where $\sigma= \sigma_n\circ\dots\circ\sigma_1$
and $\R=\{ z=0\}$, is a divisor with normal crossings.
Let put
$$\sigma^*(C)=C^{\prime}+\sum_{i=1}^n c_i E_i$$
and
$$\sigma^*(\R)=E_0+\sum_{i=1}^n a_i E_i,$$
where $C^{\prime}$, $E_0$,
and $E_i, 1\le i\le n$
are the strict transforms in $U_n$ of $C$,
$\R$, and the exceptional divisors of $\sigma_i, 1\le i\le n$, respectively.
In this notation, for
the constant $M$ mentioned above
one can take any $m$ such that
$$ m(\sum a_i E_i)-(\sum c_i E_i)$$
is a strictly positive divisor, i.e., if $ma_i-c_i > 0$ for all $i$.
Indeed, let $s$ be a singular point of $\sigma ^{-1}(C\cup \R)_{red}$,
$s\in C^{\prime}\cap E_i$ for some $i$. Choose local
coordinates $(z_i,w_i)$ in a neighbourhood of $s$ such that $z_i=0$
is an equation of $E_i$ and $w_i=0$ is an equation of $C^{\prime}$.
We have
$$\sigma^{*}(z)=z_i^{a_i}\, \, \, \, \,
\sigma^{*}(P(z,w))=z_i^{c_i}w_i$$
up to a function non-vanishing at $s$ and
$$\sigma^{*}(\overline P(z,w))=
(w_i+z_i^{Ma_i-c_i}\sum \sigma^{*}(w^{k-j}r_j(z)))z_i^{c_i}.$$
Therefore
the germ given by $\overline P(z,w)=0$ has the same resolution of
singularities as $(C,0)$ has.

{\bf 2.2. Local braid monodromy over a point (several germs).}
Now let $C\in \C^2$ be an affine reduced algebraic curve
given in coordinates $(z,w)$ by equation
$$w^m+\sum_{i=1}^{m}q_{i}(z)w^{m-i}=0, \qquad q_i\in\C[z]. \tag 2.2$$
Notice that any affine algebraic curve is given by such an
equation after a suitable linear change of coordinates.

Consider $C_{\varepsilon_1}=C\cap (D_{1}(\varepsilon_1)
\times D_{2}(\varepsilon_2))$,
$\varepsilon_2=\varepsilon_1^{-1}$.
Assuming that $0<\varepsilon_1<<1$, the projection
$\pr_{\mid C_{\varepsilon_1}} : C_{\varepsilon_1} \to D_1$,
$D_1=D_{1}(\varepsilon_1)$,
is a proper map of degree
$m$ with the unique critical
value $z=0$ (contrary
to the situation in 2.1, here
the number of critical points may be more than one).
By a traditional abuse of language we speak on
$C_{\varepsilon_1}$ as the
{\it germ of $C$ over $0$ with respect to} $\pr$.
As in the local case (see 2.1), we fix
$u\in\partial D_1$ and put
$D_{2,u}=\pr^{-1}(u)$,
$K=K(u)=\{w_1,\dots,w_m\}
=D_{2,u}\cap C$.
Giving to $\partial D_1$
its counter clock-wise orientation,
we get over $\partial D_1$
by means of $\pr_2(C\cap \partial D_1)$
an oriented loop of $m$-tuples in $D_2$ and,
thus, a braid $\tilde b_{(C_{\varepsilon_1},o)}\in
B_m=B_m[D_2,K]$.

Let $\pr_{\mid C}^{-1}(0)=\{(0,w_1^0),\dots ,(0,w_s^0) \}$. Then
the germ $C_{\varepsilon_1}$ over $0$ of $C$ splits into the
disjoint union $C_{\varepsilon_1}= \bigsqcup_i^s
C_{\varepsilon_1,i}$ of $W$-prepared germs of singularities  of
multiplicities $k_i$, $1\le i\le s$, $k_1+\dots+k_s=m$ with
centers at $(0,w_i)$, $1\le i\le s$. Let $k_1,\dots, k_t\geq 2$
and $k_{t+1}=\dots k_s=1$. We need to select a suitable polydisc
for each
of
these $W$-prepared germs.
Therefore, choose
$\varepsilon_3>0$ and adjust $\varepsilon_1<<\varepsilon_3$ so
that each $C_{\varepsilon_1,i}\subset D_1\times E_{i}$, where
$E_{i}=\{ \mid w-w_i^0\mid <\varepsilon_3 \}$ and $E_{i}\cap
E_{j}=\emptyset$ for $i\neq j$. Denote by $K_i=K_i(u)=
\{w_{i,1},\dots,w_{i,k_i}\}=D_{2,u}\cap C_{\varepsilon_1,i}$ and
$E_{i,u}=(D_1\times E_{i})\cap D_{2,u}$. The embeddings
$(E_{i},K_i)\subset (D_{2,u},K)$ defines embeddings
$\eta_i:B_{k_i}[E_{i,u},K_i]\subset B_m[D_{2,u},K]$ so that
$$ \tilde b_{(C_{\varepsilon_1},o)}=
\prod_{i=1}^t b_{(C_{\varepsilon_1,i},o)}\subset B_m.$$ As in 2.1,
let choose in each of $B_{k_i}[E_{i,u},K_i]$ as standard
generators the half twists. Then, each
$b_{(C_{\varepsilon_1,i},o) }
\in B_{k_i}$ becomes the
standard form of braid monodromy of $C_{\varepsilon_1,i}$ with
respect to $\pr$. The union of the images of these generators
under $\eta_i$ can be extended to a set of standard half-twist
generators in $B_m=B_m[D_{2,u},K]$. Thus, we get a topological
braid monodromy
$$
b_{(C_{\varepsilon_1},o)}=
b_{(C_{\varepsilon_1,1},o)}\cdot\dots\cdot b_{(C_{\varepsilon_1
,t},o)} \subset T
$$
(see the definition of $T$ in section 1.3)
and call it the {\it standard form}
of braid monodromy of $C$ over $0$.

The above construction depends only on the numbering of the points
$p_{\mid C}^{-1}(0)=\{(0,w_1^0),\dots ,(0,w_s^0) \}$ and the
extension of the union of the images of the generators under
$\eta_i$ to the sets of generators in $B_m=B_m[D_{2,u},K]$.
Therefore, the standard form of braid monodromy is
{\bf defined uniquely up to conjugation action} of $B_m$.

As it follows from the above construction,
$$\alpha (b_{(C_{\varepsilon_1},o
)})=\tilde b_{(C_{\varepsilon_1},o )}.$$

To conclude this subsection, let us recall a well known elementary
transformation replacing a germ of a pencil by a local
singularity. Namely, with a germ of $C$ over $0$ given by (2.2)
one can associate another germ $\bar C_{\varepsilon_1}$ given by
$$w^m+\sum_{i=1}^{m}z^{i}q_{i}(z)w^{m-i}=0. \tag 2.3$$
It has only one point over $0$ and
we call it the {\it associated singularity}
of $C$ over $0$.

{\bf Remark 2.1.} The geometric meaning of the associated
singularity is the following. We can consider $\C^2$ with
coordinates $(z,w)$ as an affine chart in some ruled surface $\pr:
\F_N\to \Bbb P^1$, $N>0$, non-intersecting the exceptional section
$E_N$ of $\F_N$, and $C$ as an algebraic curve contained in this
chart of $\F_N$. Perform the elementary transformation $\tau
:\F_N\to \F_{N+1}$ with center at the intersection point of $E_N$
with the fiber $\R_0$ over the point $z=0$. Then, equation (2.3)
is the equation of the image $\tau(C)\subset \F_{N+1}$ in the
corresponding chart of $\F_{N+1}$. The associated singularity
determines uniquely the germ of $C$: to get it back it is
sufficient to perform the inverse transformation $\tau^{-1}$.

The next claim
sets up a relation between
a standard braid monodromy form of a
germ of a curve
and the braid monodromy of the singularity associated to the germ.
Since this claim is not used below, we put its proof
to the end of this section.

\proclaim{Claim 2.1} Let $(C_{\varepsilon_1},o)$ be a germ of an
algebraic curve of degree $m$ over a point $o$ and $(\bar
C_{\varepsilon_1},o)$
the germ of its associated singularity. Then
$$b_{(\bar C_{\varepsilon_1},o)}=\alpha(b_{(C_{\varepsilon_1},o)})\Delta^2_m.$$
\endproclaim
{\bf Remark 2.2.} Since a germ of a curve is determined by its
associated singularity, it follows from this claim that two germs
$(C_{1,\varepsilon_1},o)$ and $(C_{2,\varepsilon_1},o)$ are
topologically equivalent if and only if
$\alpha(b_{(C_{1,\varepsilon_1},o)})=\alpha(b_{(C_{2,\varepsilon_1},o)})$.
Therefore, introducing a definition of braid monodromy
factorizations of algebraic curves (see below), we could restrict
ourselves to factorization semigroups of the first level. However,
in Section 3 we extend the class of curves under investigation
from algebraic to topological Hurwitz curves. A braid monodromy
factorization of a topological Hurwitz curve requires to consider
factorization semigroups of the second level. Therefore, for unity
of exposition, we define the braid monodromy factorization of an
algebraic curve as an element of some factorization semigroup of
the second level over a braid group as well.

{\bf 2.3. Braid monodromy factorization over a disc.} Here, as in
subsection 2.2, we consider a polynomial
$P(z,w)=w^m+\sum_{i=1}^{m}q_{i}(z)w^{m-i} \in \C [z,w]$ having no
multiple factors and the curve $C$ in $\C^2$ given by $P(z,w)=0$.

Pick any $r>0$ such that
no critical value of
$\pr\vert_C$ belongs to
the boundary of the disc $D_{1}(r)$.
Denote by $z_1, \dots ,z_n\in D_{1}(r)$
the critical values of  $\pr\vert_C$
situated in $D_{1}(r)$.
Choose a positive $\r<<1$
such that the discs $D_{1,i}(\r)=\{ z\in \C \, \mid
\,\, \mid z-z_i\mid < \r \, \}$, $i=1,\dots, n$,
would be disjoint.
Pick any points $u_i\in \partial
D_{1,i}(\r), 1\le i\le n,$ and a point
$u_0\in \partial D_{1}(r)$. Put
$D_{2,u_i}=\pr^{-1}(u_i)$, $i=0,\dots,n$, and
$K(u_i)=\{w_{i,1},\dots,w_{i,m}\} =D_{2,u_i}\cap C$.
Select disjoint simple
paths $l_i\subset \Cl D_{1}(r)\setminus \bigcup^n_1 D_{1,i}(\r)$,
$i=1,\dots,n$, starting at $u_0$ and
ending at $u_i$ and renumber the points, if needed,
in a way that the product $\gamma_1\dots \gamma_n$
of the loops $\gamma_i=l_i\circ
\partial D_{1,i}(\r)\circ l_i^{-1}$
would be equal to $\partial D_1(r)$ in $\pi_1(\Cl D_{1}(r)
\setminus \{ z_1, \dots ,z_n\}, u_0)$
(as usual, $\partial D$ states for a one counter-clock wise turn loop).

Now, as in subsection 2.2, one can associate with each $\gamma
_i$, $1\le i\le n,$ an element $b_i\in T_{\emptyset}\subset
T\subset S_{B_m}$, $B_m=B_m[D_2,K(u_0)]$, where $D_2$ is a disc of
a big radius in $w$-plane. For each $1\le i\le n$ the conjugacy
class of $b_i$ is the standard form of braid monodromy of the germ
of $C$ at $z_i$. The factorization $b_1\cdot\dots \cdot b_n\in
\Cal T$ is called a {\it braid monodromy factorization} of the
polynomial $P(z,w)$, or that of the curve $C$, over $D_{1}(r)$.
Since, for given $P(z,w)$, the braid monodromy factorizations over
$D_{1}(r)$ coincide up to Hurwitz moves and simultaneous
conjugations (see \cite{\refMT}), they all have the same
factorization type in the sense of section 1. We denote the
factorization type of braid monodromy factorizations of $P(z,w)$
over $D_1(r)$ (respectively, of the curve $C$ given by $P(z,w)=0$)
by $\text{bmt}(P(z,w), D_1(r))$ (respectively, $\text{bmt}(C,
D_1(r))$). If $P(z,w)$ has no critical values outside $D_1(r)$,
we speak
simply on {\it braid modnoromy factorization of } $P(z,w)$ (or
$C$) and denote it factorization type by
$\text{bmt}(P(z,w))$ (respectively, $\text{bmt}(C)$).

Consider all braid monodromy standard forms of all germs over $0$
of all affine algebraic curves of degree $m$. Denote by $P$ the
union of their orbits under the conjugation action of $B_m$ on
$S_{B_m}$ and put $\Cal{P}=\Cal P_m=S(S_{B_m},P)$. Note that
$A_k\subset \Cal {P}$, where $A_k$ is defined in section 1.2 (in
particular, the elements belonging to $A_0$ correspond to a simple
tangency between $C$ and $z=0$, belonging to $A_1$ correspond to a
node, and belonging to $A_2$ to an ordinary cusp). Since $\Cal
{P}=S(S_{B_m},P)$ is embedded into $\Cal T$, each braid monodromy
factorization of a polynomial $P(z,w)$ over a disc $D_{1}(r)$
defines an element in the semigroup $\Cal{P}$ and two braid
monodromy factorizations of two polynomials over discs are of the
same  braid factorization type if and only if the corresponding
elements in $\Cal{P}$ belong to the same orbit of the conjugation
action of $B_m$ on $\Cal{P}$. Therefore, we call the orbits of the
conjugation action of $B_m$ on $\Cal{P}$ {\it geometric braid
factorization types}. Note that $\beta_S :\Cal P \to S_{B_m}$ is
an embedding. Therefore, often we do not make difference between
$\Cal P$ and its image $\beta_S(\Cal P)$.

All the notions introduced in this Subsection extend literally
to any closed domain in $\C$ diffeomorphic to a closed disc.

{\bf 2.4. Polynomial realizations.} Let $\F_N$, $N\geq 1$, be a
relatively minimal ruled rational surface, $\pr: \F_N \to \Bbb
P^1$ its ruling, $\R$ a fiber of $\pr$ and $E_N$ the exceptional
section, $E_N^2=-N$. We choose a
point $\infty \in \Bbb P^1$ and put $\R_{\infty}=p^{-1}(\infty )$.

Consider the linear system
$$\Bbb P=\Bbb P_{N,m}=\Bbb PH^0(\F_N, \Cal O_{\F_N}(mE_N+mN\R)).
$$
In the chart $\F_N\setminus (\R_{\infty}\cup E_N) \simeq \C^2$
we can choose
coordinates $(z,w)$ such that the restriction of $p$
to $\F_N\setminus (\R_{\infty}\cup E_N)$ coincides with the projection
$\pr :(z,w)\mapsto z$.
With respect to these coordinates
any element $\bar C\in \Bbb P$ can be given by equation
$$P(z,w)=0,$$
where
$P(z,w)=\sum_{i=0}^{m}q_{i}(z)w^{m-i}$ and
$q_i(z)=\sum_{j=0}^{iN}a_{i,j}z^j$ are polynomials of
degrees $\leq iN$.

If $\bar C$ is an irreducible curve, then, since the intersection
number $\bar C E_N$ is zero, the
projection $\pr_{\mid \bar C}:\bar C\to \Bbb P^1$ is a proper map
of degree $m$. It has $m(m-1)N$ critical values (counted with
multiplicities) which are found from the equation
$$R(z)=0,$$
where  $R(z) =R_{P,P^\prime_w}(z) =\sum r_kz^k$ is the resultant
of $P(z,w)$ and $P^{\prime}_w(z,w)=\frac{\partial}
{\partial w}P(z,w)$, $\deg R(z)=m(m-1)N$.

The coefficients $r_k$ of the resultant admit
polynomial expressions
in coefficients of $P$. These polynomials
$r_k\in \C [a_{0,0},\dots ,a_{m,mN}]$ define a rational map
$$\Cal R:\Bbb P\to \text{Sym}_{m(m-1)N}\Bbb P^1\simeq\Bbb P^{m(m-1)N},$$
where $\text{Sym}_{m(m-1)N}\Bbb P^1$ is the symmetric product of
$m(m-1)N$ copies of $\Bbb P^1$.
The only non regular points of this map correspond to reducible
curves $\bar C\in\Bbb P$
with a multiple component different from a fiber
(here, the regularity at a point is equivalent to the existence of
a continuous extension, and the checking of the latter is
straightforward).

Below we will need the following lemmas.

\proclaim{Lemma 2.1} Any braid monodromy factorization
$b=\text{bmf}\, (\bar C) \in \Cal P$ of a generic curve $\bar C\in
\Bbb P$ is equal to $(\delta^2)^N$.
\endproclaim

\demo{Proof} By Lemma 1.2, $(\delta^2)^N$ is the unique element of
its orbit under the conjugation action of $B_m$. Thus, it is
sufficient to show that $(\delta^2)^N$ is equal to some braid
monodromy factorization. Such an equality is well-known if $N=1$
(see, for example, \cite{\refMo}). We show how to deduce the
general statement from this particular case.

Pick $z_1,\dots,z_N$ distinct from $0$.
Consider a curve $\bar C_0\in\Bbb P$ given by equation

$$\prod_{j=1}^m (wz_2\dots z_N+(-1)^Na_j(z-z_1)\dots (z-z_N))=0,
\quad a_i\neq a_j\,\,\text{if}\,\, i\neq j,$$
and a generic curve  $\bar C\in\Bbb P$
sufficiently close to $\bar C_0$ (i.e., a curve in the intersection
of a small topological neighborhood of $\bar C_0$ with
the set of curves $C\in\Bbb P$ for which
the projection $\pr_{\vert C}:C\to \Bbb P^1$ has only simple
critical points and disctinct critical values). The set of
critical values of $\pr_{\vert\bar C}$
splits into $N$ subsets $\{ z_{k,1},\dots ,z_{k,m(m-1)}\}$, $k=1,\dots ,N$,
lying in small disjoint discs
$D_{1,k}(\varepsilon)
=\{ \mid z-z_k\mid < \varepsilon \} $, $0<\varepsilon <<1$.
Correspondingly, $b= b_1\cdot\dots\cdot b_N$, where
$b_k$ is the braid monodromy factorization of $\bar C$ over
$D_{1,k}(\varepsilon)$.

Let us show that $b _1=\delta^2$ (the proof of
$b_2=\dots=b_N=\delta^2$ is the same). Define a path $\bar C_0(t),
0\le t\le 1,$ in $\Bbb P$ by equations
$$
\align
&\prod_{j=1}^m (wz_2(t)\dots z_N(t)+(-1)^Na_j(z-z_1)(z-z_2(t))\dots
(z-z_N(t)))=0,\\
&\qquad\qquad\qquad z_j(t)\to\infty\,\,\text{when}\,\, t\to 1,\quad
z_j(0)=z_j.
\endalign
$$
It connects
the curve $\bar C_0(0)=\bar C_0$ with a curve $\bar C_0(1)$
given by equation
$$\prod_{j=1}^m (w-a_j(z-z_1))=0.$$
Since $\Cal R$ is regular at each point of the path
$\bar C_0(t)$, we can find
a path $\bar C(t)$ sufficiently close to it such that:
for each $t\neq 1$ the projection $\pr:\bar C(t)\to \Bbb P^1$ has
only simple critical points
and distinct critical values, $\bar C(0)=\bar C$, $\bar C(1)$ is given
by polynomial $\sum_{i+j=m}a_{i,j}(z-z_1)^iw^{j}=0$, and for each
$0\le t<1$ all the $m(m-1)$ critical points of $\bar C(t)$ lie in
$\{ \mid z-z_1\mid < \varepsilon \} $.
It implies that the braid monodromy factorization of $\bar C$ over
$D_{1}(\varepsilon)$ coincides with
the braid monodromy factorization of $\bar C(1)$ over
$D_{1}(\varepsilon)$, which is equal to $\delta^2$
(case $N$=1).
\qed \enddemo

\proclaim{Lemma 2.2} Let $\bar C_0\in \Bbb P$
have an ordinary
singular point of multiplicity $m$ at $z=0,w=0$, i.e.,
the equation of $\bar C_0$ is of the form
$P_{\geq m}(z,w)=\sum_{i+j\geq m}a_{i,j}z^iw^{j}=0$
where the polynomial
$\sum_{i+j= m}a_{i,j}w^{j}=0$ has $m$ distinct roots. Let
$Q(z,w)\in\Bbb P$ be a polynomial such that
for some $\varepsilon_0>0$ and for all $0<\varepsilon<<1$, the curve
$\bar C_{\varepsilon}$ given by equation
$$
P_{\geq m}(z,w)+\varepsilon Q(z,w)=0
$$
has exactly $m(m-1)$ different critical values in
$D_{1}(\varepsilon_0)=\{\mid z\mid<\varepsilon_0 \}$.
Then, the braid monodromy factorization
of $\bar C_{\varepsilon}$ over
$D_{1}(\varepsilon_0)$ is equal to $\delta^2$.
\endproclaim

\demo{Proof} Connect $\bar C_0$ with a curve
$\bar C_1\in\Bbb P$ given by equation
$\sum_{i+j= m}a_{i,j}z^iw^{j}=0$ by
a path $\{\bar C_t\}_{0\le t\le 1}$ in
$\Bbb P$, where $\bar C_t$ is given by equation
$$P_t(z,w)=\sum_{i+j= m}a_{i,j}z^iw^{j}+\sum_{i+j> m}a_{i,j}(t)z^iw^{j}=0.$$
Since the leading part of the equations is nondegenerate, for any $t$ the
resultant $\Cal R(\bar C_t)$ has the following roots:
$z_{1,0}(t)=\dots =z_{m(m-1),0}(t)=0$,
$z_{k,0}(t)\neq 0$ for $k>m(m-1)$.
Hence, we can choose
$\varepsilon_0 >0$ such that
$2\varepsilon_0 < \mid z_{k,0}(t)\mid $ for $k>m(m-1)$ and for
all $t\in [0,1]$. For $\varepsilon$
small enough the image
$\{\Cal R(\bar C_{t,\varepsilon})\}$ of the path
$\{P_{t}(z,w)+\varepsilon Q(z,w)=0\}$
is contained in a neighborhood $V_{\varepsilon_0}\subset
\text{Sym}_{m(m-1)N}\Bbb P^1$ of $\{\Cal R(\bar C_{t})\}$,
$$
\align
V_{\varepsilon_0}=\{ (z_1,\dots ,z_{Nm(m-1)}) \,\, \mid &
\, \, \mid z_i \mid <\varepsilon_0\, \,
\text{for}\, \, i\leq m(m-1) \\
 &  \text{and} \, \,\mid z_i\mid > \varepsilon_0 \, \,
 \text{for}\, \, i> m(m-1) \} .
\endalign
$$
Changing slightly the coefficients of the polynomial $Q(z,w)$,
we can assume that
for $0<\varepsilon<<1$ and any $0\le t\le 1$,
$\Cal R(\bar C_{t,\varepsilon})$ belongs to
$V_{\varepsilon_0}\cap \{ z_i\neq z_j \, \,
\text{for}\, \,
i,j\leq m(m-1), \, \, i\neq j \}$.
Under this assumption, all the curves
$\bar C_{t,\varepsilon}$ have the same braid monodromy factorization
over $D_{\varepsilon_0}$.

The Viro patch-working, see \cite{\refV}, is based on
quasi-homogeneous changes of coordinates exclusively, and thus
respects the braid monodromy. Hence, to finish the proof it
remains to replace $\bar C_{t,\varepsilon}$ by a generic
polynomial obtained by patch-working generic polynomials
$$\sum_{i+j< m}a'_{i,j}z^iw^{j}+\sum_{i+j= m}a_{i,j}z^iw^{j},\quad 
\sum_{i+j= m}a_{i,j}z^iw^{j}+\sum_{i+j> m}a'_{i,j}z^iw^{j}$$ close
to $\bar C_1$ and to apply Lemma 2.1 to the first one (case
$N=1$). \qed \enddemo

\proclaim{Lemma 2.3}
Let $\bar b$ be
the braid monodromy factorization of
an affine part $C$
of a projective curve $\overline C\in \Bbb P_{N_1,m}$
over a disc $D_{1}(r)$, $r>>1$.
Then, there is $M=M_C\in\N$ such that for
any $N_2\geq M$ there is a projective
curve $\widetilde C_1\in \Bbb P_{N_1+N_2,m}$ such that
$\widetilde b=\bar b\cdot (\delta^2)^{N_2}$ is
the braid monodromy factorization of its affine part $C_1$
over the disc $D_{1}(r_1)$, $r_1>>r$ and $\bar b$ is
the braid monodromy factorization of $C_1$
over the disc $D_{1}(r)$.
\endproclaim

\demo{Proof}
Denote by $z_1,\dots, z_k$
the critical values of $\pr_{\mid C}$ and
by $(z_i,w_{i,1}),\dots, (z_k,w_{i,t_{i}})$
the critical points over $z=z_i$.
Attribute to the germ $(C_{i,j},(z_i,w_{i,j}))$
of $C$ at $(z_i,w_{i,j})$
a number $M_{i,j}$ defined in section 2.1. Put
$M_i =\max_j M_{i,j}$ and $M_C=\sum_i M_i$.
Choose points $z_{k+i}$, $i=1,\dots, N_2$
such that $\mid z_{k+i}\mid >>r$ and perform
$N_2> M_C$ elementary transformations
with centers at the intersection points
of the fibers $z=z_{k+i}, 1\le i\le N_2,$
with the exceptional sections of $\F_{N_1+i-1}$,
$i=1,\dots, N_2$.
Denote by $\widetilde C\subset \F_{N_1+N_2}$
the strict transform of $\bar C$ under
the composition $\tau :\F_{N_1}\to \F_{N_1+N_2}$
of these transformations and by
$$\widetilde P(z,w)=w^m+\sum_{i=1}^{m}q_{i}(z)w^{m-i}=0$$
the equation of its affine part.
The only critical values of
$\pr_{\mid \widetilde C}$ are
$$
z_1,\dots, z_k, z_{k+1}, \dots , z_{k+N_2}.
$$
Over $z=z_{k+i}$
there is only one point of $\widetilde C$
and at this point
the curve $\widetilde C$ is equivalent, with respect to $\pr$, to
$$ \prod_{i=1}^m(w-c_i(z-z_{k+i}))=0.$$
The braid monodromy factorization over $D_{1}(r)$ of $\widetilde
C$ coincides with the braid monodromy factorization of $C$, since
the centers of the transformation $\tau$ lie over $z_{k+i}$, $1\le
i\le N_2$, and $\mid z_{k+i}\mid >>r$. By the choice of the
constants $M_j$ and by Lemma 2.2, a curve
$C_{\overline{\varepsilon}}$ given by
$$\widetilde P(z,w)+ \prod_{i=1}^k(z-z_i)^{M_i}
(\sum_{j=0}^{m-1}\varepsilon_jw^j) =0,$$
where all $\varepsilon_j$ are
generic and close to zero,  has
the braid monodromy factorization
$\widetilde b=\bar b\cdot (\delta^2)^{N_2}$
over the disc $D_{1}(r_1)$, $r_1>>r$. Indeed, by Bertini's theorem,
the curves $C_{\overline{\varepsilon}}$ are non-singular over
$\Bbb P^1\setminus \{z_1,\dots, z_k\}$ for almost all
$\overline{\varepsilon}$. Moreover, if for some $\overline{\varepsilon}$
a fiber $\R$ over a point lying near some $z_{k+i}$
is tangent the curve $C_{\overline{\varepsilon}}$
with multiplicity greater than 2 at some point,
then we can choose $\overline{\varepsilon}_1$
close to $\overline{\varepsilon}$ such that $C_{\overline{\varepsilon}_1}$
and $\R$ have the only one simplest tangent point. Therefore
there is $\varepsilon_0$ such that
for almost all $\overline{\varepsilon}$ sufficiently close to zero,
exactly $m(m-1)$ distinct critical values of $\pr:
C_{\overline{\varepsilon}}\to \Bbb P^1$ lie in each neighborhood
$\{\mid z-z_{k+i}\mid <\varepsilon_0\}$.
\qed\enddemo

\proclaim{Lemma 2.4}
If $b\in \Cal{P}$ is the braid monodromy factorization
over a disc $D$ of an affine curve
$P(z,w)=w^m+\sum_{i=1}^{m}p_{i}(z)w^{m-i}=0$
with $p_i(z)=\sum_{j=0}^{iN}a_{i,j}z^j, N\ge 1,$ then there exists
an affine curve
$Q(z,w)=w^m+\sum_{i=1}^{m}q_{i}(z)w^{m-i}=0$
such that
\roster
\item"{($i$)}" $q_i(z)=\sum_{j=0}^{iN}b_{i,j}z^j$;
\item"{($ii$)}" the polynomial $w^m+\sum_{i=1}^{m}b_{i,iN}w^{m-i}$
has $m$ different roots;
\item"{($iii$)}" all critical points of $Q(z,w)$ lying over
the complement $\C\setminus D$  are non-degenerate, in particular,
the curve $C$
given by $Q(z,w)=0$ is non-singular over $\C\setminus D$;
\item"{($iv$)}" a braid monodromy
factorization of $C$ over the disc $D$
is equal to $b$.
\endroster
\endproclaim

\demo{Proof} The second property, (ii), can be achieved by a generic choice
of the fiber at infinity. After that, it is sufficient to apply to
$$Q(z,w)+ \prod_{i=1}^k(z-z_i)^{M_i}
(\sum_{j=0}^{m-1}\varepsilon_jw^j) =0,$$
the same
inductive correction procedure as at the end of the proof of Lemma 2.3
\qed
\enddemo

\proclaim{Theorem 2.1} For any
$b\in \Cal{P}$ there are $N\in \N$ and a polynomial
$P(z,w)=w^m+\sum_{i=1}^{m}q_{i}(z)w^{m-i}$ such that
\roster
\item"{($i$)}" $q_i(z)=\sum_{j=0}^{iN}a_{i,j}z^j$;
\item"{($ii$)}" the polynomial $w^m+\sum_{i=1}^{m}a_{i,iN}w^{m-i}$
has $m$ different roots;
\item"{($iii$)}" all critical points of $P(z,w)$ lying over
the complement $\C\setminus D_{1}(1)$ of the disc
$D_{1}(1)=\{ z\in \C \, \mid
\,\, \mid z\mid < 1 \, \}$ are non-degenerate, in particular,
the curve $C$
given by $P(z,w)=0$ is non-singular over $\C\setminus D_{1}(1)$;
\item"{($iv$)}" $b$ is a braid monodromy
factorization of $P(z,w)$ over the disc $D_{1}(1)$.
 \endroster
\endproclaim

\demo{Proof} Let
$b=\rho (g_1^{-1})(b_1)\cdot\dots\cdot
\rho (g_n^{-1})(b_n)$, where
$b_j\in T$ are standard forms of braid
monodromies of algebraic germs $C_j$ over points
$z=z_j$, $\mid z_j\mid <1$. Let
$$w^m+\sum_{i=1}^{m}q_{j,i}(z-z_j)w^{m-i}=0$$
be an equation of the germ $C_j$ at $z=z_j$ and
$$q_j(z,w)=w^m+\sum_{i=1}^{m}(z-z_j)^{i}q_{j,i}(z-z_j)w^{m-i}=0 $$
be the equation of the singularity $\bar C_j$ associated with $C_j$.

\proclaim{Proposition 2.1}
There is a polynomial $Q(z,w)=w^m+\sum_{i=1}^{m}\widetilde q_{i}(z)w^{m-i}$
having at each $(z_j,0)$ the same type of singularity
as $\bar C_j$ and satisfying  conditions \rom ($i$\rom)--\rom ($iii$\rom) of
Theorem 2.1 for some $N$.
\endproclaim

\demo{Proof} Choose a constant $M>M_{(\bar C_j,z_j)}$ for all $j=1,\dots, n$,
where $M_{(\bar C_j,z_j)}$ is the constant defined in section 2.1.

In the beginning we construct a polynomial $\overline Q(z,w)$ of degree
$\deg_w \overline Q=m$ having singularities at $(z_j,0)$, $j=1,\dots, n$,
of the same
types as $(\bar C_j,z_j)$. The polynomial $\overline Q(z,w)$ will be
constructed by $n$ steps. Put $\overline Q_1(z,w)=q_1(z,w)$.
Assume that we have constructed a
polynomial $\overline Q_k(z,w)$ of degree
$\deg_w \overline Q_k=m$ having singularities at
$(z_j,0)$, $j=1,\dots, k$, of the same
types as $(\bar C_j,z_j)$.
Consider a polynomial $q_{z_1,\dots z_k}(z)=
((z-z_1)\dots (z-z_k))^M$. We have  $q_{z_1,\dots z_k}(z_{k+1})\neq 0$.
Therefore we can find polynomials $p_{i}(z)$ such that
the polynomial
$$\overline Q_{k+1}(z,w)=\overline Q_k(z,w)+
 \sum_{i=1}^{m}q_{z_1,\dots z_k}(z)p_{i}(z)w^{m-i}$$
 have the same types  of singularities at
$(z_j,0)$, $j=1,\dots, k+1$,
as $(\bar C_j,z_j)$. Indeed,
$\overline Q_{k+1}(z,w)$ has the same singularities
as $(\bar C_j,z_j)$ for $j=1,\dots, k$,
whatever is the choice of the polynomials $p_i$,
since $M$ is big enough, and choosing appropriate
coefficients for the polynomials $p_{i}(z)$, the polynomial
$\overline Q_{k+1}(z,w)$ will have the same type of singularity at
$(z_{k+1},0)$ as $(\bar C_{k+1},z_{k+1})$.

To complete the prove of Proposition 2.1,
it remains to apply Lemma 2.4.
\qed\enddemo

By Remark 2.1, it follows from Proposition 2.1 that there is a
polynomial $Q(z,w)$ satisfying conditions ($i$)--($iii$) of
Theorem 2.1 for some $N=N_1$ and such that the set of critical
values of $\pr\vert_{Q(z,w)=0}$ lying in $D_{1}(1)$ coincides with
$\{ z_1,\dots, z_n \}$ and over $z=z_j$, $j=1,\dots, n$, the germ
of the curve $C$ given by $Q(z,w)=0$ has the same topological type
as $(C_j,z_j)$. Therefore, the braid monodromy factorization of
$Q(z,w)$ over the disc $D_{1}(r)$, $r>>1$, is  equal to
$$\bar b=\rho (\bar g_1^{-1})(b_1)\cdot\dots\cdot
\rho (\bar g_n^{-1})(b_n)\cdot s,$$ where $s\in S_{A_0}$. Besides,
from conditions ($i$)--($iii$) of Theorem 2.1 it follows that
$C\subset \overline C$ is an affine part of a projective curve
$\overline C\in \F_{N_1}$, such that all singular points of
$\overline C$ belong to $C$.

It follows from Theorem 1.1, Proposition 2.1 and Lemma 2.3 that
there are $N_2>>1$ and a polynomial $\overline P(z,w)$ of
$\deg_w=m$ whose braid monodromy factorization over the disc
$D_{1}(r)$, $r>>1$, is equal to $b\cdot s\cdot (\delta^2)^{N_2}$,
where $s$ is an element belonging to $S_{A_0}$. Therefore there is
a domain $U \subset\C$ diffeomorphic to the disk $D_{1}(1)$ such
that the braid monodromy factorization  of $\overline P(z,w)$ over
$U$ is equal to $b$. There is an analytic isomorphism of $U$ and
$D_{1}(1)$. Therefore, there is a polynomial $\widetilde
P(z,w)=w^m+\sum_{i=1}^{m}q_{i}(z)w^{m-i}$ with holomorphic in
$D_{1}(1)$ coefficients $q_{i}(z)$ such that the braid monodromy
factorization  of $\widetilde P(z,w)$ over $D_{1}(1)$ is equal to
$b$. Now, it remains to approximate the holomorphic functions
$q_i(z)$ by polynomials with the same jets
(of sufficiently big degree) at the critical values of $p$
belonging to $D_1 (1)$ and then to apply Lemma 2.4. \qed\enddemo

{\bf Remark 2.3.} For given $b\in \Cal P$, such that $\beta (b)=
\Delta^2$, let
$$P(z,w)=w^m+\sum_{i=1}^{m}q_{i}(z)w^{m-i}$$
be a polynomial whose existence is stated in Theorem 2.1.
As it follows from Proposition 3.1,
there does not exist
an upper bound for degrees of the polynomials $q_i(z,w)$
depending only on the topological types of critical points of $P(z,w)=0$
with respect to $\pr$
even in the case of the types $A_n, \,n\le 2$
(simple tangents, nodes and ordinary cusps). On the other hand,
as it follows from the proof of Theorem 2.1, there is an effective
bound on $\deg_z\bar P$
for the intermediate polynomial $\bar P$, which provides a given braid
monodromy factorization over some domain.

{\bf Remark 2.4.} The case $N=1$ covers the case of algebraic
curves in $\C P^2$. In the latter case to define a braid monodromy
factorization one choose a point outside the curve, and the
blowing up reduces this case to the case of curves in $\F_1$.
\newline
{\it Proof of Claim 2.1.} Let realize a germ
$(C_{\varepsilon_1},o)$ of an algebraic curve over a point $o$ as
a germ of some projective curve $\widetilde C\subset \F_N$.
Let $b(\widetilde C)$ be 
a braid monodromy factorization of $\widetilde C$. 
Then $b(\widetilde
C)=b_{(C_{\varepsilon_1},o)}\cdot\widetilde b$, where $\widetilde
b\in \Cal P$ is the product of the braid monodromies over all the
critical values of the projection $\pr$ except $o$. We have
$$\beta(b(\widetilde C))=\beta(b_{(C_{\varepsilon_1},o)})\beta(\widetilde
b)=\Delta_m^{2N}.$$ Perform the elementary transformation $\tau$
with center at $\R_o\cap E_{\infty}$. The curves
$\tau^{-1}(\widetilde C)$ and $\widetilde C$ have the same
critical values, the braid monodromy factorization of
$\tau^{-1}(\widetilde C)$ is $b(\tau^{-1}(\widetilde C))=b_{(\bar
C_{\varepsilon_1},o)}\cdot\widetilde b$, and also
$$\beta(b(\tau^{-1}(\widetilde C)))=\beta(b_{(\bar C_{\varepsilon_1},o)})\beta(\widetilde
b)=\Delta_m^{2(N+1)}.$$ Therefore, $\beta(b_{(\bar
C_{\varepsilon_1},o)})=\beta(b_{(C_{\varepsilon_1},o)})\Delta_m^{2}$.
It remains to note that $\beta(b_{(\bar
C_{\varepsilon_1},o)})=b_{(\bar C_{\varepsilon_1},o)}$ and
$\beta(b_{(C_{\varepsilon_1},o)})=\alpha
(b_{(C_{\varepsilon_1},o)})$ in view of the identifications
made
in section 1.3. \qed

\subheading{\S 3.\ Braid monodromy factorizations of Hurwitz
curves}

{\bf 3.1. Hurwitz curves.} As in section 2,
let $\F_N$ be a relatively minimal ruled rational surface,
$N\geq 1$, $\pr: \F_N \to \Bbb P^1$
the ruling, $\R$
a fiber of $\pr$ and $E_N$ the exceptional section, $E_N^2=-N$.

\proclaim{Definition 3.1} The image $\bar H=f(\Si)\subset \F_N$ of
a smooth map $f:\Si \to \F_N\setminus E_N$ of an oriented closed
real surface $\Si$ is called a Hurwitz curve \rom(in $\F_N$\rom)
of degree $m$ if there is a finite subset $Z\subset\bar H$ such
that: \roster \item"{($i$)}" $f$ is an embedding  of
the surface $\Si\setminus f^{-1}(Z)$ and for any
$s\notin Z$, $\bar H$ and the fiber $\R_{\pr(s)}$ of $\pr$ meet at
$s$ transversely and with positive intersection number;
\item"{($ii$)}" for each $s\in Z$ there is a neighborhood
$U\subset \F_N$ of $s$ such that $\bar H\cap U$ is a complex
analytic curve, and the complex orientation of $\bar H\cap U
\setminus \{s\}$ coincides with the orientation transported from
$\Si$ by $f$; \item"{($iii$)}"  the restriction of $\pr$ to $\bar
H$ is a finite map of degree $m$.
\endroster
\endproclaim

For any Hurwitz curve $\bar H$ there is one
and only one minimal $Z\subset\bar H$ satisfying the conditions
from Definition 3.1. We denote it by $Z(\bar H)$.

A Hurwitz curve $\bar H$ is called {\it cuspidal } if for
each $s\in Z(\bar H)$ there is a
neighborhood $U$ of $s$ and local analytic coordinates $x,y$ in
$U$ such that
\roster
\item"{($iv$)}" $\pr_{\mid U}$ is given by
$(x,y)\mapsto x$;
\item"{($v$)}"  $\bar H\cap U$ is given by
$y^2=x^k, k\ge 1$.
\endroster
It is called {\it ordinary cuspidal} if $k\leq 3$ in ($v$) for
all $s\in Z(\bar H)$, and {\it nodal}
if $k\leq 2$.

Since $\bar H\cap E_N=\emptyset$, one can define a braid monodromy
factorization $b(\bar H)\in {\Cal P}$ of $\bar H$ as in the
algebraic case. For doing this, we fix a fiber $\R_{\infty}$
meeting transversely $\bar H$ and consider $\bar H\cap \C^2$,
where $\C^2=\F_N\setminus (E_N\cup \R_{\infty})$. Choose $r>>1$
such that $\pr (Z)\subset D_{1}(r)\subset\C= \C P^1\setminus\pr
F_\infty$, $Z=Z(\bar H)$. Denote by $z_1,\dots,z_n$ the elements
of $\pr (Z)$. Pick $\r$, $0<\r<<1$, such that the  discs
$D_{1,i}(\r)=\{ z\in \C \, \mid \,\, \mid z-z_i\mid < \r \, \}$,
$i=1,\dots, n$, would be disjoint. Select arbitrary points $u_i\in
\partial D_{1,i}(\r)$ and a point $u_0\in \partial D_{1}(r)$. Put
$D_{2,u_i}=\pr^{-1}(u_i)$, $i=0,\dots,n$, and
$K(u_i)=\{w_{i,1},\dots,w_{i,m}\}=D_{2,u_i}\cap H$. Choose
disjoint simple paths $l_i\subset \Cl D_{1}(r)\setminus
\bigcup^n_1 D_{1,i}(\r)$, $i=1,\dots,n$, starting at $u_0$ and
ending at $u_i$ and renumber the points in a way that the product
$\gamma_1\dots\gamma_n$ of the loops $\gamma_i=l_i\circ \partial
D_{1,i}(\r)\circ l_i^{-1}$ would be equal to $\partial D_1(r)$ in
$\pi_1(\Cl D_{1}(r)\setminus \{ z_1, \dots ,z_n\}, u_0)$.

As in Section 2, each
$\gamma _i$ defines an element $b_i\in T\subset S_{B_m}$. The
factorization $b_1\cdot\dots \cdot b_n\in \Cal T$ is called a {\it
braid monodromy factorization} of $\bar H$.
In fact,
each $b_i$ is conjugated to a braid monodromy standard form of
some algebraic germ over $z_i$. Hence, $b(\bar H)=b_1\cdot
\dots\cdot b_n$ belongs to ${\Cal P}$ (see Section 2.3 for the
definition of $\Cal P$). The orbit of this element under the
conjugation action of $B_m$ in $\Cal P$ is called the {\it
geometric braid factorization type} and denoted by bmt.

\proclaim{Lemma 3.1} For a Hurwitz curve $\bar H\subset \F_N$
it holds
$$\beta (b(\bar H))=\Delta^{2N}.$$
\endproclaim

\demo {Proof} Over $\Bbb P^1\setminus D_{1}(r)$ the curve $\bar H$
is the union of $m=\deg \bar H$
pairwise disjoint sections of $\pr$.
By an isotopy of sections they can be transformed into the sections defined,
in affine coordinates, by $w=v_iz^N, 1\le i\le m, v_i\in\C.$
The latter sections form the braid $\Delta^{2N}$ over $\partial D_1(r)$.
\qed\enddemo

The converse statement is also proved in a straightforward way.

\proclaim{Theorem 3.1} \rom(\cite{\refMoi}\rom) For any
$b=b_1\cdot \dots\cdot b_n\in {\Cal P}$ such that $\beta
(b)=\Delta^{2N}$ there is a Hurwitz curve $\bar H\subset \F_N$
with a braid monodromy factorization $b(\bar H)$ equal to $b$.
\endproclaim

{\bf 3.2. Hurwitz isotopies.}
The following definition generalizes
the notion of Hurwitz curves. This generalization corresponds to
replacement of $\Cal P$ by $\Cal T$.

\proclaim{Definition 3.2}
The image $\bar H=f(\Si)\subset \F_N$ of a
continuous map $f:\Si \to \F_N\setminus E_N$
of a smooth oriented closed real surface $\Si$
is called a topological Hurwitz curve \rom(in $\F_N$\rom) of degree $m$
if there is a finite subset $Z\subset\bar H$ such that:
\roster
\item"{($i$)}" $f$ is a smooth embedding
of the surface $\Si\setminus f^{-1}(Z)$ and for any $s\notin Z$,
$\bar H$ and the fiber $\R_{\pr(s)}$ of $\pr$ meet at $s$
transversely and with positive intersection number;
\item"{($ii$)}"  the restriction of $\pr$ to $\bar H$ is a finite map
of degree $m$. \rom {(We call a map finite if the preimage of each
point is finite.)}
\endroster
\endproclaim

As in the case of Hurwitz curves, there is one and only one
minimal $Z$, which we denote by $Z(\bar H)$. We say that $\bar H$
is {\it $Z$-generic}
(with respect to $\pr$) if the
restriction of $\pr $ to $Z(\bar H)$ is injective.

\proclaim{Definition 3.3} Two Hurwitz \rom(respectively,
topological Hurwitz\rom) curves $\bar H_1$ and $\bar H_2\subset
\F_N$ are called {\it $H$-isotopic} if there is a fiberwise
continuous isotopy $\phi_t:\F_N\to \F_N$, $t\in [0,1]$,  smooth
outside the
the fibers $F_{\pr(s)}$, $s\in Z(\bar H_1)$, and %
such that
\roster
\item"{($i$)}" $\phi_0=\text {Id}$;
\item"{($ii$)}" $\phi_t(\bar H_1)$ is a Hurwitz \rom(respectively,
topological Hurwitz\rom) curve for all $t\in [0,1]$; \item"{($iii$)}"
$\phi_1(\bar H_1)=\bar H_2$; \item"{($iv$)}" $\phi_{t}(E_N)= E_N$
for all $t\in [0,1]$.
\endroster
\endproclaim

If $\bar H_1$ is a topological Hurwitz curve, at any $p\in Z(\bar
H_1)$ there is a well-defined ($W$-prepared) germ $(D,H_1=\bar
H_1\cap D,\pr)$ of this curve in 
a bi-disc
$D= D_1\times D_2$, $D_1=D_1(\epsilon_1)$, $D_2=
D_2(\epsilon_2)$, $0<\epsilon_1<<\epsilon_2$, centered at $p$ and
such that the restriction of $\pr$ to $H_1$ is a proper map of
finite degree. If $\epsilon_1,\epsilon_2$ are sufficiently small,
then: $F_{\pr(p)}\cap H_1=p$;
the above degree does not depend on
$\epsilon_1,\epsilon_2$;
and, in the same way as in the algebraic
case, the link $\partial D\cap H_1$ defines a unique, up to
conjugation, braid $b\in B_k$, where $k$ is the above degree. So
that, we may speak on a {\it $tH$-singularity} $(D,H_1,\pr)$  of
{\it degree} $k$, of {\it type} $b$ and of {\it interlacing
number} $l=l(b)$.

When we are given a link $L\subset\partial D_1\times D_2$
realizing a braid $b\in B_k$, we associate with it a {\it standard
conical model} of a topological singularity of type $b$. It is
given by $H=C(L)$,
$$C(L)=\{ (rz,rw)\,\, \vert \,0\le r\le 1, (z,w)\in L \}. $$

As is well-known, if ($D,C$) is a germ of a
$W$-prepared analytic singularity then the germ
($D,C,\pr$) is homeomorpic to the cone
singularity of type $b=\pr^{-1}(\partial D_1\cap C$.

It is convenient to describe a $tH$-singularity $(D,H_1,\pr)$ by
means of a flow. To make the corresponding formulae more
transparent let put $\epsilon_1=\epsilon_2=1$ and move $H_1$ by a
suitable (homothety like) fiber preserving $H$-isotopy into the
cone body $|w|<\rho |z|$, $\rho=1-\epsilon$. Consider a smooth map
$v$ from $D_1\setminus\{0\}$ to the space of smooth vector fields
on $D_2$ such that $v(z)(w)=-w$ for any $|w|\ge\rho |z|$. Then the
flow defined by
$$
\frac{dz}{dt}=-z,\quad \frac{dw}{dt}=v(z)(w)
$$
transforms any braid $L\subset \partial D_1(1)\times D_2(\rho)$
into a $tH$-singularity. Reciprocally, any $tH$-singularity
$(D,H_1,\pr)$ can be represented in such a way as soon as $H_1$ is
contained in $|w|<\rho |z|$. For example, the standard conical
model $C(L)$ is given by any pair $(v,\rho)$ like above with
$v(z)(w)=-w$ for $|w|<\rho |z|$.

\proclaim{Claim 3.1}
For any $tH$-singularity $(D=D_1\times
D_2,H_1,\pr)$ there is an $H$-isotopy \rom(preserving each fiber\rom)
identical over $\partial D$ and transforming the
singularity in its standard conical model $C(L)$, $L=H_1\cap
(\partial D_1\times D_2)$.
\endproclaim

\demo{Proof} Let represent the $tH$-singularity $(D=D_1\times
D_2,H_1,\pr)$ and its conical model $C(L)$ by two flows as above:
one flow is associated with a pair $(v_0,\rho)$ and another with $(v_1,\rho)$.
Consider the family of $tH$-singularities given by the
pairs
$(v_t=tv_1+(1-t)v_0,\rho)$ and denote by $\Phi_{t,z}$ the imbeddings
$\{\frac{z}{|z|}\}\times D_2\to \{z\}\times D_2$ given by the flow
associated with $v_t$. To accompany the constructed family of
$tH$-singularities by some ambient $H$-isotopy $\phi_t$ connecting
$H_1$ with $C(L)$ over $D_1$, it is sufficient
to take the flow associted with the vertical vector fields
$\frac{d}{dt}(\Phi_{t,z}\circ\Phi_{1,z}^{-1})$ extending them for
each $t$ by $0$ to the whole $D_1\times D_2$. \qed\enddemo

{\bf Remark 3.1.} In fact, the procedure we used in the proof of
Claim 3.1 gives, as well, an $H$-isotopy between any two
$tH$-singularities with the same link $L$. If the both
singularities are singularities of Hurwitz curves, they remain
singularities of Hurwitz curves (i.e., algebraic singularities)
during the isotopy. Certainly, this isotopy is not necessary
smooth at the singular point.

The definitions of braid monodromy factorizations and  braid
factorization types extend literally  from Hurwitz to topological
Hurwitz curves.

Note only that if $(D,H,\pr)$ is a $tH$-singularity,
then due to Claim 3.1 it is determined uniquely, up to $H$-isotopy
identical on $\partial D$, by its boundary closed braid
$b=H\cap \pr^{-1}(\partial D_1)$ and respectively by its standard
tbmf-form (see section 1.3). Therefore the braid monodromy of a
germ of a topological Hurwitz curve over a point is naturally
defined as an element of $T$, and the braid monodromy
factorizations of topological Hurwitz curves as elements of $\Cal
T$. They satisfy the relation $\beta(b(\bar H))= \Delta^{2N}$. The
inverse statement, that each factorization $b\in\Cal T$ with this
property is realized by a topological Hurwitz curve, has the same
proof as Theorem 3.1.

The braid group $B_k=B[D_2,\{w_1,\dots,w_k\}]$, with
$w_1,\dots,w_k\in D_2$ acts on $\pi_1 = \pi_1(D_2\setminus \{
w_1,\dots ,w_k\} )$ in a natural way. We say that the action of
$b\in B_k$ on $\pi_1$ is {\it inseparable} if only the elements of
the subgroup of $\pi_1$ generated by $\partial D_2$ are fixed
under the action of $b$. Standard generators $a_1,\dots, a_{k-1}$
of $B_k$ being fixed, we mean by a {\it geometric base} of
$\pi_1(D_2\setminus \{ w_1,\dots ,w_k\} )$  any set of generators
$\{ x_1,\dots, x_k\} $ in which the natural action of $B_k$ is
given by standard formulas, i.e., $a_i(x_j)=x_j$ for $j\ne i, i+1$
and $a_i(x_{i+1})=x_{i}$, $a_i(x_{i})=x_ix_{i+1}x_i^{-1}$ for
$i\le k-1$. For any such base, $\partial D_2 =x_1\dots x_k$ .

\proclaim{Lemma 3.2} Let $b\in B_k$ be an inseparable element and
$\{ x_1,\dots, x_m\} $ a geometric base of $\pi_1(D_2\setminus \{
w_1,\dots ,w_m\} )$. Regard $b$ as an element of
$$B_{k,0}=B[D_2,\{w_1,\dots,w_k\}]\subset
B[D_2,\{w_1,\dots,w_m\}]=B_m$$
and consider the induced action of $b$ on
$\pi_1 = \pi_1(D_2\setminus \{ w_1,\dots ,w_m\} )$.
Then, the subgroup $F(b)$ of the fixed
elements in $\pi_1$ under the action of $b$ is generated by
$l=x_1\dots x_k$ and $x_{k+1},\dots, x_m$.
\endproclaim

\demo{Proof}
Evidently,
$l$ and $x_{k+i}$ with any $i\ge 1$ belong to
$F(b)$.
Let $g\in F(b)$. Write $g$ as a reduced word in letters
$\{x_1,\dots ,x_m\}$ and their inverses
$$ g=s_1s_2\dots s_n,$$
where $s_i$ are reduced words which are
non-empty if $i\ne 1$
and which are words in $\{x_1, \dots x_k\}$ and their inverses
if $i$ is odd, and
in $\{x_{k+1}, \dots x_m\}$ and their inverses
if $i$ is even.
Since such a representation is unique, we deduce from $g\in F(b)$
that each $s_i$
belongs to $F(b)$. Now the result follows from the definition
of inseparable elements.
\qed\enddemo

A $tH$-singularity $(D,H,\pr)$ is said to have an {\it inseparable type}
if its type $b\in B_k$ is inseparable. Note that
a $tH$-singularity of inseparable type $b$
has a tbmf-form $(b,\bold 1_{\emptyset})$.

\proclaim{Lemma 3.3}
The singularities
of algebraic curve given by $w^k=z^n$
are inseparable for all $n\geq 1$
and $k\geq 1$.
\endproclaim

\demo{Proof}
The type of such singularity equals to $b=(a_1\dots
a_{k-1})^n\in B_k$, where $a_1,\dots , a_{k-1}$ is a set of
standard generators of $B_k$. The equality $\Delta_k^2=(a_1\dots
a_{k-1})^k$ implies $b^k=\Delta_k^{2n}$. On the other hand, the
action of $\Delta_k^{2}$ coincides with the conjugation action by
the element $\partial D_2\in \pi_1$. Now it remains to note that
the centralizer of any $g\ne 1$ in a free group coincides with the
maximal infinite cyclic subgroup containing $g$.
\qed\enddemo

One can show that Claim 3.1 can not be extended literarly to any
germ of a topological Hurwits curve over a point. In fact, it does
not hold at least if the germ consists of several connected
components one of which is not of inseparable type. The following
example is a simplest one.
\newline
{\bf Example 3.1.} Put $m=3$ and consider
$H_1, H_2$ such that:
$bmf(H_1)=bmf(H_2)=(\bold 1,\bold 1_{\{1,2\}})$;
$H_1$ and $H_2$
are rotation symmetric (i.e., invariant under $(z,w)\mapsto
(e^{i\phi}z,w)$);
in $H_1$ one, figure-V, component is obtained by
rotation of a cone over $2$ points (lying in $\{
\text{Im}\, z=0\}\times D_2$; this component corresponds to $\bold
1_{\{1,2\}}$) and the other one by rotation of a string not linked
with the above cone, while in $H_2$ the latter string is linked
with one, and only one, branch of the cone.

\proclaim{Theorem 3.2}
$Z$-generic topological Hurwitz curves with
singularities of inseparable types
are $H$-isotopic if and only if
they have the same braid monodromy type.
\endproclaim
\demo{Proof} The necessity part is obvious.

Assume that $Z$-generic topological Hurwitz curves $\bar H_1, \bar
H_2$ have singularities of inseparable type and that $bmt(\bar
H_1)=bmt(\bar H_2)$. The latter implies that there is an
$H$-isotopy which transforms $\bar H_1$ into $\bar H_2$ everywhere
except over a union of disjoint disc neighborhoods of the points
of $\pr Z(\bar H_1)$, which are transformed by this isotopy in a
union of disjoint disc neighborhoods of $\pr Z(\bar H_2)$, cf.
\cite{\refKT}. It remains to solve the problem of extension of an
$H$-isotopy from a boundary of a disc inside the disc. Here, the
assumption on the type of singularities is essential.

So, assume that the traces $H_1, H_2$ of $\bar H_1,\bar H_2$ over
a disc $D_1=D_1(\epsilon)$ are contained in $D=D_1\times
D_2(r)$, coincide over the boundary of $D_1$,
that their braid monodromy $b\in T$ over $\partial D_1$ is of
inseparable type and that $D_1$ contains only one point of $Z(H_1)$ and
only one
point of $Z(H_2)$. Without loss of generality, we may assume that
$Z(H_1)=Z(H_2)=\{(0,0)\}$.
Now, it is sufficient to construct an $H$-isotopy identical over
boundary which transforms $H_1, H_2$ at least over one interval
connecting the boundary of $D_1$ with $0$, say over the radius
$I_-=\{\text{Re}\, z\le 0, \text{Im}\, z=0\}\subset D_1$, and to
extend it over a neighborhood of $0$. After that the $H$-isotopy
can be easily extended to the remaining part of $D_1$, since it is
a topological disc over which there are no singular points.

Let $l\ge 2$ be the interlacing number
equal to the
degree of the singularity,
since $\bar H_i$ have only
the singularities of inseparable type.
Due to the definition of $T$,
there is a standard tbmf-form $(\widetilde  b_1, \bold
1_{\emptyset})$,
$b_1\in B_{l,0}$, conjugated to $b$,
and due to the definition of the braid monodromy factorization,
$H_1$ and $ H_2$ split each into
$k+1=m-l+1$ connected components $H_{i,1},\dots,H_{i,k+1}$,
$i=1,2$, and these splittings coincide over $\partial
D_1$. Each component is a topological cone; $H_{i,1}$ is
a cone over the braid $b_1\in B_{l,0}$ and the other are cones
over $1$-braids $b_{2}\in B_{1,l},\dots,b_{k+1}\in B_{1,m-1}$.
By Claim 3.1, we may assume that $H_{1,1}=H_{2,1}=C(L)$
where $L$ is the trace of $ H_{1,1}=H_{2,1}$ in $\partial
D_1\times D_2(r)$.
In addition, we may assume that the
trace of $C(L)$ over some $D_1(\epsilon')$, $\epsilon^{\prime}<<\epsilon$,
is contained in
$D_1(\epsilon')\times \Int D_2(r^{\prime})$, $r^{\prime}<<r$, and  that
$H_{1,2},\dots,H_{1,k+1}$ (respectively,
$H_{2,2},\dots,H_{2,k+1}$) are distinct constant sections $w=w_j$,
$l+1\le j\le m$ (respectively, $w=w_j'$) of $\pr$ over
$D_1(\epsilon')$ which are contained in $D_1(\epsilon')\times A_2$,
$A_2=\{r^{\prime}\le \vert w\vert \le r\}\subset D_2(r)$. By an
$H$-isotopy the values $w_j$ (respectively, $w_j'$) of these
constant sections can be arbitrary continuously changed following
arbitrary braid not going through $D_2(r')$; we call such
a change {\it twisting-untwisting}.
Note that, in particular, when
$w_j=w_j'$ for all $j\ge l+1$ we achieve $H_1=H_2$ over a
neighborhood of $0$. Thus, it is sufficient to show that there
exists an $H$-isotopy identical on the boundary which transforms
$H_1$ in $H_2$ over $I_-$.

To describe
the action of $b_1$ on $\pi_1(D_2(r)\setminus \{w_1,\dots,w_l\},
w_0)$, $w_0\in\partial D_2(r)$, let us replace the standard
projection $\pr_2 : T\to D_2(r)$, $T=A_1\times D_2(r)$,
$A_1=\{\epsilon'\le \vert z\vert\le \epsilon\}\subset D_1$, by a
fibration $\pr':T\to D_2(r)$ 
which
coincides with $\pr_2$ on $\partial D_1(\epsilon')\times
D_2(r)\cup A_1\times\partial D_2(r)$, with the standard projection
$C(L)\to L$ on $C(L)\cap T$
given by cone structure, and
which is constant on each of $H_{1,j}, j\ge 2$. Such a fibration
is given by the flow defined by any vector field on $T$
which
is tangent to $C(L)$, all $H_{1,j}$,
$j\ge 2$, and $A_1\times\partial D_2(r)$ and obtained as a lift of
a rotation invariant radial nowhere zero vector field on $A_1$.
Then $\pr'(C_-)$, $C_-= C(L)\cap \pr_1^{-1} I_-$, consists of $l$
points, which we denote by $w_1,\dots,w_l$. Let $W$ be the element
of $\pi_1(T\setminus C(L), (-\epsilon', w_0))$ given by the
constant section $w=w_0$ over $\partial D_1(\epsilon')$. If
$x=\pr'_*(y)$ where $ y\in\pi_1(T\setminus C(L), (-\epsilon',
w_0))$ is realized by a loop lying over $I_-$, the image $x^b$ of
$x$ under the action of $b_1$ is determined by $x^b=\pr'_*(y^b)$,
$yW=Wy^b$.

Now, consider a loop $\bar y_j, l+1\le j\le m,$ which starts at
$(-\epsilon', w_0)$ descends to $(-\epsilon', w_j)$ without
entering into $D_2(r')$, goes along $H_{1,j-l+1}$ over
$I_-$, returns along $H_{2,j-l+1}$ and goes up to $(-\epsilon',
w_0)$ without entering into $D_2(r')$. Rotating $I_-$ we
get a map of torus for which $\bar y_j$ is a meridian-map, and $W$
is given by a parallel-map. Hence, $y_j$ realized by $\bar y_j$
and $W$ commute, so that $x_j=\pr'_*y_j$ is invariant under the
action of $b_1$. Since $b_1$ is inseparable,
for any $j$ there is an isotopy of $D_2(r)\setminus \{ w_1,\dots,
w_{l}\}$ identical on the boundary which
transforms $\pr'\bar y_j$ in a loop outside $D_2(r')$.

Finally, to transform all $H_{1,j}$ in $H_{2,j}$ over $I_-$ by an
$H$-isotopy which is identical on $\partial D \cup C(L)$ let apply
the following induction procedure. At the first step by an
$H$-isotopy identical on $\partial D$ and $C(L)$ we transform
$\pr'\bar y_2$ in a loop outside $D_2(r')$, and then after
defined by this loop twisting-untwisting transform it in a
constant loop, i.e., identify $H_{2,2}$ with $H_{1,2}$ over $I_-$.
When $H_{1,j}=H_{2,j}$ over $I_-$ for any $j\le i$ ($i<m-l+1$), we
consider the element realized in
$\pi_1(D_2(r)\setminus\{w_1,\dots,w_l,\dots,w_{l+i-1}\})$ by
$\pr'\bar y_{i+1}$. Since, by the same argument as above, it is
also invariant under the natural action of $b_1$. Therefore
by Lemma 3.2, it is a product
of elements corresponding to $\partial D_2(\epsilon')$ and loops
around $w_j$, $j>l$, not entering into $D_2(\epsilon')$. Thus, again
as above, we can apply an $H$-isotopy realizing the homotopy of
$\pr'\bar y_{i+1}$ in a loop outside $D_2(r')$ and then
after twisting-untwisting an $H$-isotopy transforming it in a
constant loop.
\qed
\enddemo

\proclaim{Corollary 3.1} Two Hurwitz curves $\bar H_1, \bar H_2
\subset \F_N$ with singularities of inseparable types are isotopic if
$\beta_S(\text{bmt}(\bar H_1))=\beta_S(\text{bmt}(\bar H_2))$.
\endproclaim

\demo{Proof} By isotopy, we can slightly move the points belonging
to $Z(\bar H)$ to get a $Z$-generic Hurwitz curve having
 $\text{bmt}=\beta_S(\text{bmt}(\bar H))$.
Now the proof follows from Theorem 3.2.  \qed\enddemo

\proclaim{Corollary 3.2} $Z$-generic Hurwitz curves with
singularities of types $w^k=z^n$ are $H$-isotopic if and only if
they have the same braid monodromy type.\qed
\endproclaim

\demo{Proof} It follows from Lemma 3.3.
\qed\enddemo

As is shown in  \cite{\refKT}, in the case of $Z$-generic cuspidal
Hurwitz curves the $H$-isotopy can be made smooth. The theorem
below is a generalization of this result to the non generic case.

\proclaim{Theorem 3.3} Let $\bar H_1$ and $\bar H_2\subset \F_N$
be cuspidal Hurwitz curves. Then
$\bar H_1$ and $\bar H_2$ are smoothly $H$-isotopic if and only if
$\bar H_1$ and $\bar H_2$ have the same braid monodromy type.
\endproclaim

\demo{Outline of proof}
The braid monodromy standard form of a germ
of a cuspidal curve over a point having singularities of types
$w^2=z^{n_i}$ is defined by the collection of exponents
$\{n_1,\dots,n_t\}$ and it is of the form
$$(a_1^{n_1},\bold 1_{\emptyset})\cdot (a_3^{n_2},\bold
1_{\emptyset})\cdot \dots \cdot (a_{2t-1}^{n_t},\bold
1_{\emptyset}),$$ where $2t\leq m$.

It is not hard to prove (we omit the
proof) that the centralizer $C(b)\subset B_m$ of
$$b=a_1^{n_1}a_3^{n_2}\dots a_{2t-1}^{n_t}$$
is generated  by $a_1,a_3,\dots ,a_{2t-1}$,
$a_{2t+1},a_{2t+2},\dots, a_{m-1}$ and some elements
$c_{i}$, $1\leq i\leq t$,
and $d_{i,l}$, $i\ne l$, $1\leq i,l\leq t$, where
$$c_{i}=a_{2t}\dots a_{2i+1}(a_{2i}a_{2i-1}^2a_{2i})a_{2i+1}^{-1}\dots
a_{2t}^{-1}$$ and
$$d_{i,l}=(a_{2l}a_{2l-1})\dots
(a_{2i-2}a_{2i-3})
(a_{2i}a_{2i-1}a_{2i+1}a_{2i})^2
(a_{2i-2}a_{2i-3})^{-1}\dots (a_{2l}a_{2l-1})^{-1}$$
if $l<i$ and
$$d_{i,l}=(a_{2l-2}a_{2l-1})\dots
(a_{2i+2}a_{2i+1})
(a_{2i}a_{2i+1}a_{2i-1}a_{2i})^2
(a_{2i}a_{2i+1})^{-1}\dots (a_{2l-2}a_{2l-1})^{-1}$$
if $l>i$.

Now as in the proof of Theorem 3.2, by an $H$-isotopy, we identify
$\bar H_1$ and $\bar H_2$ over the complement of the union of
small discs around the points of $\pr(Z(\bar H_1))$.
In each such disc
$D_1(\epsilon)$ pick a much more small disc $D_1(\epsilon^{\prime})$
over which we can identify $H_1$ and $H_2$
with the same disjoint union of several sections and several singularities
given by $(w-w_i)^2=z^{n_i}$.
Select a trivialization of $\pr$ over
$D_1(\epsilon)\setminus D_1(\epsilon^{\prime})$
with respect to which
the trace $H_1$ of $\bar H_1$ is
radialy constant. Recall that
$H_2$ and $H_1$ coincide over
the boundary of $D_1(\epsilon)\setminus D_1(\epsilon^{\prime})$.
Since $H_1$ and $H_2$ have the same $tbmf$-form
$$\widetilde b=(a_1^{n_1},\bold 1_{\emptyset})\cdot (a_3^{n_2},\bold
1_{\emptyset})\cdot \dots \cdot (a_{2t-1}^{n_t},\bold
1_{\emptyset}),$$ the braid $g=\hat r\cap H_2 \in B_m$, where
$\hat r= \{ z\in D_1(\epsilon)\, \mid \, z\in \Bbb R,\,
\epsilon^{\prime}\leq  z\leq \epsilon\}$, commutes with
$$b=a_1^{n_1}a_3^{n_2}\dots a_{2t-1}^{n_t}.$$

According to the
description of the centralizer $C(b)$ of $b$, the element $g\in
C(b)$ can be written as a word $\hat g$ in the alphabet consisting
of the generators of $C(b)$ and their inverses. The first letter of this word
can be removed by one of the twisting-untwisting
isotopies: any $c_{i}$ and $a_l$ with
$l\geq 2t+1$ (and their inverse) can be removed
by twisting-untwisting isotopy used in the proof of Theorem 3.2,
$a_{2i-1}$ with $1\leq i\leq t$ can be removed by
twisting-untwisting isotopy used in \cite{\refKT}, and every
$d_{i,l}$ can be removed by a ``composition'' of these two
isotopies. After a sequence of such isotopies the connected
components of $H_2$ will not be linked with the connected
components of $H_1$, i.e., when we will have $g=\bold 1$. Now the
end of the proof coincides with the corresponding part of the
proof of Theorem 3.2. \qed\enddemo

{\bf 3.3. Almost-algebraic curves.} By Theorem 2.1,
for any element $b \in {\Cal P}$ there are $N\in \N$ and
an algebraic curve $\bar C\subset \F_N$
not intersecting $E_N$ (and intersecting transversely $\R_\infty$)
and such that its
affine part $C\subset \F_N\setminus (E_N\cup \R_\infty)$
has over $D_{1}(1)$
the braid monodromy factorization equal to $b$.
For $b \in {\Cal P}$, the minimal
number $N=N(b)$ such that $b$ is a braid monodromy factorization
of an algebraic curve
$\bar C\subset \F_N$ over
$D_{1}(r)\subset \Bbb P^1$
with some $r>0$
is called the {\it realizing degree} of $b$.
If $\beta (b)=\Delta^{2k}$, the number
$$d(b)=N(b)-k$$
is called the {\it deficiency} of $b$.

\proclaim{Definition 3.4}
A Hurwitz curve $\bar H\subset \F_N$ is called
an almost-algebraic curve
if $\bar H$ coincides with
an algebraic curve $C$
over $D_{1}(r)$ and
with the union of $m$
pairwise disjoint sections
$H_{\infty, 1},\dots ,H_{\infty, m}$ of $\pr$
over $\Bbb P^1\setminus D_{1}(r)$.
\endproclaim

The following  theorem is a direct
corollary of
Theorems 2.1, 3.2, and 3.3.

\proclaim{Theorem 3.4} \roster
\item"{($i$)}" For each $b\in \Cal P$ with $\beta
(b)=\Delta^{2N}$ there is an almost algebraic Hurwitz curve $\bar
H\subset \F_N$ whose braid monodromy factorization $b(\bar H)$ is
equal to $b$.
\item"{($ii$)}" Any $Z$-generic
Hurwitz curve $\bar H\subset \F_N$ with singularities of types
$w^k=z^n$ \rom(in general, with singularities of inseparable
type\rom) is $H$-isotopic to an almost-algebraic curve. If the
Hurwitz curve is cuspidal, this isotopy can be chosen smooth.
\endroster \qed
\endproclaim

A braid monodromy factorization
$b(\bar H)=b_1\cdot \dots\cdot b_n\in {\Cal P}$
of a generic ordinary cuspidal Hurwitz curve $\bar H\subset \F_N$
has the following form:
$$b(\bar H)=\prod_{i=1}^n(q_ia_1^{r_i}q_i^{-1})\in {\Cal P},$$
where $r_i$ depends on the type of singularity $s_i\in Z$:
$r_i = 3$ if $s_i$ is a cusp
(i.e., $\bar H$ is given locally at $s_i$ by
$y_i^2=x_i^3$),
$r_i = 2$ if $s_i$ is a node
(i.e., given by
$y_i^2=x_i^2$), and $r_i = 1$ if $s_i$ is a tangency point
(i.e., given by $y_i^2=x_i$).

\proclaim{Proposition 3.1} There is an infinite sequence
$\bar H_i\subset \F_1$, $i\in \N$,
of generic
ordinary
cuspidal almost-algebraic curves of degree $54$
with exactly 378 cusps and 756 nodes
and of pairwise distinct braid monodromy types.
In particular, they are not $H$-isotopic
to each other, almost all of them are not
isotopic to an algebraic cuspidal curve and, moreover,
$$\lim _{i\to \infty} d(b(\bar H_i))=\infty .$$
\endproclaim

\demo{Proof} The existence of a sequence of $\bar H_i$ like in the
first part of the statement follows from Theorems 2.1 and 3.4
applied to the series of pairwise distinct cuspidal braid
factorization types of $\Delta^2_{54}$ found by Moishezon, see
\cite{\refMoi} Theorem 1.

Let fix $N$ and introduce a space $B$
of algebraic curves
$$
\bar C \in
\Bbb P=\Bbb PH^0(\F_N, \Cal O_{\F_N}(54E_N+54N\R))
$$
which do not contain $E_N$,
have no critical points over the boundary of the disc
and have: $54(54-1)$ critical values lying in
$D_{1}(1)\subset \Bbb P^1$ (counting with multiplicities);
among them 378 critical values corresponding to the
cusps, 756 values corresponding to the nodes
and all the other corresponding to simple tangencies.
To prove that $\lim d(b(\bar H_i))=\infty $,
it is sufficient to show that $B$
consists of finite number of
components.
Note that $B$ is contained in
the space $\Cal M$ of all curves
in $\Bbb P$
having at least 378 cusps and 756 nodes, which is
a quasi-projective variety.

As in the proof of Lemma 2.2, consider the map
$\Cal R:\Bbb P\to \text{Sym}_{54(54-1)N}\Bbb P^1\simeq\Bbb P^{54(54-1)N}$.
It is well-defined at each point of $B$ and
$\Cal R(B)$
is contained in
$$
\align
V=\{ (z_1,\dots ,z_{54(54-1)N}) \,\, \mid &
\, \, \mid z_i \mid <1\, \,
\text{for}\, \,
i\leq 54(54-1) \\
 &  \text{and} \, \,\mid z_i\mid > 1 \, \,
 \text{for}\, \, i> 54(54-1) \} ,
\endalign
$$
so that $B\subset \Cal R^{-1}(V)\cap\Cal M$.
The set $\Cal R^{-1}(V)\cap\Cal M$ is semi-algebraic.
Let consider the semi-algebraic stratification of $\Cal M$
according equisingularity. Intersecting it with
$\Cal R^{-1}(V)$ we get a finite (semi-algebraic) stratification
of $\Cal R^{-1}(V)$. The intersection of $B$ with each stratum
of the latter
stratification is open and closed in the stratum, hence $B$ has only finite
number of connected components.
\qed\enddemo

\proclaim{Problem}
Whether any nodal almost-algebraic curve $\bar H \in \F_N$
has the deficiency $d(b(\bar H))=0$, i.e.,
$H$-isotopic
to an algebraic nodal curve $\bar C\in \F_N$?
\endproclaim

{\bf 3.4. Several remarks.}
Let
$$b(\bar H)=\prod_{i=1}^n \lambda (q_i)(b_{i})\in {\Cal P}$$
be a braid monodromy factorization of a Hurwitz curve
$\bar H\subset \F_N$, $\deg \bar H=m$,
where each $b_i$ is a standard form of braid
monodromy of an algebraic curve over a point.
Every $b_i$  is factorized into a product
$$ b_i=b_{i,1}\cdot \dots \cdot b_{i,t_i}, $$
where:
$b_{i,1}\in B^+_{k_{i,1},0}\subset B_m$,
$B^+_{k_{i,1},0}$ stating for the semigroup generated by
the elements $a_1,\dots , a_{k_{i,1}-1}$; $b_{i,2}\in B^+_{k_{i,2},k_{i,1}}$,
$B^+_{k_{i,2},k_{i,1}}$ being generated by $a_{k_{i,1}+1}$, $\dots$,
$a_{k_{i,1}+k_{i,2}-1}$; \dots ; $b_{i,t_i}\in
B^+_{k_{i,t_i},k_{i,1}+\dots +k_{i,t_{i}-1}}$,
$B^+_{k_{i,t_i},k_{i,1}+\dots +k_{i,t_{i}-1}}$ being generated by 
$a_{k_{i,1}+\dots +k_{i,t_i-1}+1}$, $\dots$ , $a_{k_{i,1}+\dots
+k_{i,t_i}-1}$. Recall that $k_{i,j}\geq 2$ for $1\leq j \leq
t_i$.

By means of this factorization one can describe the homotopy type
of the complement $\F_N\setminus (\bar H\cup E_N\cup \R_{\infty})$
(cf., \cite{\refL}). Namely, it has the homotopy type of the two
dimensional Cayley complex of the following (van Kampen - Zariski)
presentation of the fundamental group
$$\pi_1(\F_N\setminus (\bar H\cup E_N\cup \R_{\infty}))
=\langle x_1,\dots, x_m : R_{i,j,k}\rangle;
$$
in terms of the canonical action of $B_m$
on the free group with  generators $x_1,\dots, x_m$, the
defining relations $R_{i,j,k}$
($1\le i\le n$, $1\le j\le t_i$,
$k_{i,1}+\dots +k_{i,j-1}+1\leq k\leq k_{i,1}+\dots +k_{i,j}-1$)
take the form
$$q_i^{-1}(b_{i,j}(x_k))=q_i^{-1}(x_k).$$

Notice also that using the arguments from \cite{\refACC} one can
show that the braid monodromy type of topological Hurwitz curves
in $\F_N$ is preserved by some special ambient homeomorphisms.
Namely, let associate with each topological Hurwitz curve  $\bar
H$ a stratified union $\widetilde H=\bar H \cup (\bigcup_{s\in Z}
\R_{\pr(s)})\cup \R_{\infty} \subset \F_N$ (and call it an {\it
equipped curve}); two topological Hurwitz curves $\bar H_1, \bar
H_2$ has the same braid monodromy type if (and only if) there is a
homeomorphism of oriented pairs $h:(\F_N,\widetilde H_1) \to
(\F_N,\widetilde H_2)$ which transforms $\bar H_1$ in $\bar H_2$
and keeps invariant $\R_\infty$.

The condition that $h$ preserves orientation on the components of
topological Hurwitz curves is essential. This follows from the
examples (\cite{\refKK}) of irreducible (ordinary) cuspidal
algebraic curves $C_1, C_2\subset \F_1$ having different braid
monodromy types, but for which there exists a diffeomorphism
$c:(\F_1,\widetilde C_1)\to (\F_1,\widetilde C_2)$ not preserving
orientation on the components of the equipped curves (in these
examples, $c$ is the complex conjugation).

\subheading{\S 4.\ Braid monodromy factorizations of symplectic
surfaces}

In this section we consider singular symplectic surfaces in $\C
P^2$. As is known, the symplectic structure on $\C P^2$ is unique
up to symplectomorphisms and multiplication by a constant factor.
On the other hand, up to our knowledge, its uniqueness up to
isotopy and constant factor is an open question. Moreover,
rescaling of the symplectic structure is an important ingredient
of our proof of Proposition 4.2 below. These are the reasons why,
in what follows, we speak on isotopy classes of symplectic
structures.

We treat only the case of surfaces with isolated singularities and
define a {\it singular symplectic surface} $C\subset \C P^2$ as
a triple $(C,J,\omega)$ such that: $C$ is
the image $C=f(\Si)$ of an
almost everywhere
injective 
$J$-holomoprhic map
$f: \Si\to\C P^2$
of a closed Riemann surface $\Si$; 
and 
$J$ is an almost complex structure
defined on a neighborhood of $C$ and tamed by a given symplectic
structure $\omega$ of $\C P^2$.
Two singular symplectic
surfaces $(C_0=f_0(\Si),J_0^{loc},\omega_0), (C_1=f_1(\Si),J_1^{loc},\omega_1)$ are said {\it weakly
symplectically (smoothly) isotopic} if $(f_0,J_0^{loc},\omega_0),
(f_1,J_1^{loc},\omega_1)$ can be included in a continuous
(respectively, smooth) family of tamed almost complex and
symplectic structures $J_t^{loc}, \omega_t$ and that of
$J_t^{loc}$-holomorphic 
maps
$f_t:\Si\times [0,1]\to \C P^2$
such that the maps $\phi _{t} :C_0\to C_t$ given by $f_0(s)\mapsto
f_t(s)$, $s\in \Si$, are well-defined homeomorphisms for all $t\in
[0,1]$. The same definitions are applied to symplectic surfaces in
any manifold. When the ambient symplectic structure is not
changing one speaks on {\it symplectic isotopy}. In particular,
two singular symplectic surfaces are symplectically isotopic, if
there is a symplectic diffeotopy transforming one into another.

{\bf Remark 4.1.} In the case of arbitrary symplectic manifolds it
is natural to expect that the classes of symplectically isotopic
surfaces and the classes of weakly symplectically isotopic
surfaces are not the same in general. By contrary, in the case of
symplectic structures on $\C P^2$ these notions coincide. Indeed,
in this case any weakly symplectic isotopy can be rescaled
$\omega_t\mapsto\omega_t'=\lambda_t\omega_t,\lambda_t\in\Bbb R_+,$
(without changing $J_t$) to an isotopy with $\omega_t'$ not
changing their cohomology class; then it remains to apply the
Mozer theorem to get a diffeotopy making $\omega_t'$ constant; the
same diffeotopy is then applied to $f_t$ and $J_t^{loc}$.
Therefore, when we need to prove that some symplectic surfaces in
$\C P^2$ are symplectically isotopic, it is sufficient for us to
check that they are weakly symplectically isotopic.
It is this strategy that we adopt below.

{\bf Remark 4.2.} The space of almost complex structures on a
finite dimensional real vector space which are tamed by a given
symplectic form is contractible, see \cite{\refGr}. It implies
that any $\omega$-tamed almost complex structure $J^{loc}$ defined
in an open subset $U$ of a symplectic manifold $(V,\omega)$ can be
extended from a smaller neighborhood $U_0\subset U$ to a
$\omega$-tamed almost complex structure $J$ on the whole $V$.
Moreover, by the same reason, if $U=U_t, J^{loc}=J^{loc}_t,
\omega=\omega_t$ depend smoothly on one or several parameters $t$
the extensions from $U_0=U_{0,t}$ to $\C P^2$ can be chosen
depending smoothly on $t$. Therefore, one gets the same notions if
in the above definitions the locally defined $J$-structures are
replaced by tamed almost structures defined on the whole $\C P^2$.
The choice we made is motivated by simplifications in some of the
proofs. Note also that the above extension properties remain true
if, in addition, we will restrict ourself to almost complex
structures for which a given symplectic surface is
$J$-holomoprhic. In particular, in the definition of weakly
symplectically isotopic surfaces it is sufficient to have
$J_t^{loc}$ to be defined only near the singular points of $C_t$.
\vskip0.07in

As is known (see \cite{\refSi}, \cite{\refMW}), the singularities
of pseudoholomorphic curves in an almost complex four-manifold are
equivalent to the singularities of genuine complex curves in
complex surfaces up to $C^1$ coordinate change.

If a singular symplectic surface $C$ is given together with a
symplectic structure $\omega$ and a $\omega$-tamed almost complex
structure $J^{loc}$ on a neighborhood $U$ of $C$ in $\C P^2$, we
extend $J^{loc}$ from a smaller neighborhood $U_0\subset U$ to a
$\omega$-tamed almost complex structure $J$ on $\C P^2$, then
consider a generic pencil $L$ of $J$-lines, see \cite{\refGr} and
\cite{\refSik}, and define the associated braid factorization type
$bmt(C,\omega, J,L)$ as in section 2.3 in the algebraic case.

Note that  any $b\in\Cal P$ such that $\alpha(b)=\Delta^2$ (and
only the such ones) can be realized as a braid monodromy
factorization of a symplectic surface with respect to a (non
necessarily generic) pencil. Such a realization is obtained from
an almost-algebraic curve given by Theorem 3.4
by rescaling the
standard pencil where the curve is situated.

\proclaim{Proposition 4.1} The braid factorization type
$bmt\,(C,\omega,J,L)$ with respect to a gene\-ric $L$ depends only
on the symplectic isotopy class
of $C$. \rom{( In particular, it does not depend on the extension
$J$ of $J^{loc}$ to the whole $\C P^2$.)}
\endproclaim

\demo{Proof} Let $C_0, C_1$ be weakly symplectically isotopic
symplectic surfaces equipped with $\omega_i (i=0,1)$-tamed
$J$-structures, $J_0^{loc}, J_1^{loc}$, their extensions $J_0,
J_1$, and generic pencils, $L_0, L_1$, as in the definitions of
singular symplectic surface and its braid monodromy factorization,
respectively. Pick a weakly symplectic isotopy $(C_t, \omega_t,
J_t^{loc})$ and extend $J_t^{loc}$ to an isotopy $J_t$ on the
whole $\C P^2$. For each $t$, the space $(\C P^{2}_t)^*$ of
$J_t$-lines is diffeomorphic to $\C P^2$ and is equipped with the
canonical dual $E$-structure, see \cite{\refSik}. The lines
tangent to $C_t$ form the dual $E$-curve $C^*_t$. It has a finite
number of singular points (which correspond to double or excess
tangents) and the dual to $C^*_t$ is $C_t$. Hence, the nongeneric
pencils correspond to a choice of the center of the pencil $L_t$
belonging neither to $C_t$ nor to the finite number of $J_t$-lines
which are dual to the singular points of $C^*_t$. Therefore, there
is a path $L_t$ of generic pencils connecting $L_0$ with $L_1$.
Clearly, the braid monodromy factorization $C_t$ defined by $L_t$
is not depending on $t$ and the result follows. \qed
\enddemo

Recall that the braid monodromy factorization type of $m$ complex
lines in $\C P^2$ in general position coincides with $\widetilde
\delta_m^2$ defined in section 1.4 (see, for example,
\cite{\refMT}).

\proclaim{Corollary 4.1} The braid monodromy factorization type of
a generic family of $m$ pseudo-lines is equal to $\widetilde
\delta_m^2$.
\endproclaim

\demo{Proof} According to Barraud \cite{\refBa}, any generic
family of $m$ pseudo-lines is symplectically isotopic to $m$ true
lines in general position.\qed
\enddemo

\proclaim{Corollary 4.2}
There are infinitely many ordinary cuspidal symplectic plane curves of the
same degree with the same number of cusps and the same number of
nodes but two-by-two not symplectically isotopic.
\endproclaim

\demo{Proof} It follows from Propositions 3.1 and 4.1.
\qed\enddemo

{\bf Remark 4.3.} The statement reverse to Proposition 4.1 holds
at least if the symplectic isotopy is replaced by a topological
one and if, in addition, it is assumed that all singularities are
of inseparable types. In deed, as it follows from \cite{\refLM},
pencils determined by two almost complex structures tamed by a
same symplectic structure or, more generally, by symplectic
structures in a continuous family, are isotopic. Together with
Theorem 3.2 it implies that, if braid factorization types of two
singular symplectic (with respect to a same structure or with
respect to structures from the same connected component) surfaces
are equal and the singularities of the surfaces are of inseparable
types, then the surfaces are topologically isotopic. There is some
confusion in replacing a topological isotopy by a smooth one. This
is because in the smooth category a smooth isotopy is supposed to
be an ambient one, contrary to our choice in the definition of
smooth symplectic isotopies. Certainly, the ambient isotopy can be
made smooth outside the singular points, and by contrary, in
general, can not be made smooth at the singular points, since they
can have moduli even with respect to smooth changes of
coordinates.

{\bf Remark 4.4.} In the above Corollary, the pseudo-lines can be
replaced by Hurwitz curves (recall that according to our
definitions, Hurwitz curves are situated in a pencil of ordinary
lines) realizing the generator of $H_2(\C P^2)$, since after
rescaling a transversal pencil, the Hurwitz surfaces becomes
$J$-curves, so becomes pseudo-lines. As Schevchishin communicated
to us, Barraud result can be generalized to nodal symplectic
curves (without negative nodes) of genus $\le 3$. Then, the above
arguments are applied to such curves as well.

The following proposition is a partial inverse of Proposition 4.1.

\proclaim{Proposition 4.2} Two symplectic, with respect to the
Fubini-Studi symplectic structure $\omega_0$, ordinary cuspidal
surfaces are symplectically $C^1$-smoothly isotopic in $\C P^2$ if
and only if they have the same braid
factorization type with respect to a generic pencil.
\endproclaim

\demo{Proof}
The
necessity part, i.e.,
coincidence of braid monodromy factorization types, follows from
Proposition 4.1.

According to Remark 4.1, it is sufficient to find a weakly
symplectic isotopy.

Let $C_0$ be an ordinary cuspidal symplectic surface, which is a
$J_0$-holomorphic curve where $J_0$ is an almost complex structure
on $\C P^2$ compatible with $\omega_0$. By a continuous variation
of $J_0$ with a support in a neighborhood of some point $p$, make
$J_0$ integrable in a smaller neighborhood $U'_0$. Then, consider
a generic pencil $L_0$ of pseudo-lines with a center $p_0\in
U'_0$. Choose local coordinates $x,y$ near the critical points of
$C_0$ (by {\it critical} points we mean the singular points of
$C_0$ and the points of tangency between $C_0$ and $L_0$) so that
locally the elements of the $J_0$-pencil are given by the fibers
of $(x,y)\mapsto x$ and $C_0$ is defined by equations $y^2=x^k$
with $k=1,2,3$, respectively to the cases of tangency points,
nodes, and cusps; in all the cases the direction $y=0$ can be
chosen symplectic. It allows us to replace $J_0$ by a
$\omega_0$-tamed almost complex structure $J'_0$ with respect to
which the above local coordinates become $J$-holomorphic (so that
$J'_0$ is integrable near the critical points and near $p_0$) and
$C_0$, as well as the ruling $L_0$, remain $J$-holomoprhic. Since
$J_0$ and $J'_0$ are $\omega_0$-tamed and since $C_0$ is
$J$-holomorphic with respect to both of them, they can be joined
by a homotopy almost complex structures keeping $C_0$ to be
$J$-holomorphic. By Proposition 4.1, $\text{bmt}(C_0)$ does not
change. The next weakly symplectic isotopy consists in a
continuous variation of $\omega_0$, it has a support in $U'_0$ and
replaces  the pair $(\omega_0, J'_0)$ by a pair $(\omega'_0,
J'_0)$ which is standard (i.e., K\"ahler flat) in $U_0\subset
U'_0$. Such a variation is given in \cite{\refMP}, Lemma 5.5B.

On the other hand, by Theorem 3.1, there exists a Hurwitz curve
$C_1\subset\C P^2$ whose braid monodromy type (with respect to a
generic pencil $L_1$ of ordinary lines) $\text{bmt}(C_1)$ is equal
to $\text{bmt}(C_0)$ (recall that the latter is defined by means
of a generic pencil of $J_0$-lines). By rescaling a generic pencil
$L_1$ of ordinary lines and moving its center to $p_0$, we can
assume that $C_1$ is  a $J_1$-holomorphic curve, where $J_1$ is a
suitable almost complex structure tamed by $\omega_0$ and
identical with the standard one in a neighborhood of $p_0$ and in
some neighborhoods of the critical points of $C_1$.

In addition, as before, we can replace
$(\omega_0,J_1)$ by a pair $(\omega'_1, J'_1), J'_1=J_1,$ standard near
$p_0$. To proof Proposition 4.2, it is sufficient to show that
$C_0$ and $C_1$ are weakly symplectically $C^1$-smoothly isotopic
in $\C P^2$.

Then, proceed as in the proof of Proposition 4.1: pick a path
$(\omega'_t, J'_t)$ where the almost complex structures $J'_t$ are
tamed by $\omega'_t$ and integrable near $p_0$, and consider a
family of $J'_t$-holomorphic pencils of $J'_t$-lines connecting
$(J'_i,L_i)$, $i=0,1$.
Make, by a continuous
variation, the almost complex structures $J'_t$ integrable near
the critical points. Now, it remains to construct an isotopy
between $C'_0$ and
$C_1$ and to enchance it
so that it becomes a weakly symplectic isotopy. We construct the
isotopy in two steps.

At the first step,  by a diffeotopy of the pencils holomorphic
near $p_0$ and near the critical points we get a diffeotopy $C'_t$
of the curve $C_0=C'_0$. It provides us with surfaces $C'_t$ which
are holomorphic near the critical points and gives a kind
of $H$-isotopy: outside the critical points each surface $C'_t$
meets the pseudo-lines transversely and with positive intersection
number; the projection of $C'_t$ to the base of its pencil is a
finite ramified covering; near the critical points it is a complex
analytic curve which in local analytic coordinates $x,y$ such that
the projection is given by $(x,y)\mapsto x$ is defined by equation
$y^2=x^k$ with $k=1,2,3$. The resulting curve $C'_1$ is a genuine
Hurwitz curve, since $L_1$ is a pencil of ordinary lines.
Obviously, $\text{bmt}(C'_1)=\text{bmt}(C_0)$.

In the second step, the arguments from \cite{\refKT} can be
applied to the two surfaces $C'_1$ and $C_1$ having the same braid
monodromy factorization type with respect to $L_1$ to construct an
$H$-isotopic family of Hurwitz curves $H_t$, connecting $H_0=C'_1$
and $H_1=C_1$, with all the properties enumerated above.

To produce a weakly symplectic isotopy by means of the combined
isotopy ($C'_t$ followed by $H_t$) constructed above, let apply to
$(C'_t,\omega'_t)$ followed by $(H_t,\omega'_1)$ the following
rescaling of $\omega'_t$. Like at the beginning of the proof, use
the construction from \cite{\refMP} to make, by a continuous
variation of $\omega'_t$, the pairs $(\omega'_t,J'_t)$ coinciding
with the standard flat pair $(\Omega, i)$ in a small ball
$B(\delta), \delta>0$, around $p_0$. Then, consider $\Bbb
C^1$-fibrations $h_t : \C P^2\setminus\{p_0\}\to S^2$ whose fibers
$h_t^{-1}(v),v\in S^2,$ are symplectic and which coincide with the
$J'_t$-rulings $\pr: \C P^2\setminus\{p_0\}\to S^2$ outside
$B(\delta)$ and with the ordinary lines ruling in a smaller ball
$B(\delta')$. As soon as such fibrations are given, it remains to
replace $\omega'_t$: outside a smaller ball $B(\delta')$ by
$\omega'_t+N h_t^*\omega_{S^2}$, where $N\ge 0$ is a sufficiently
big constant, $\omega_{S^2}$ is a volume form on $S^2$; and inside
$B(\delta')$ by some $\Omega_N$ with $\Omega_N=\omega'_t+N
h_t^*\omega_{S^2}= \Omega+Nh^*\omega_{S^2}$ near $\partial
B(\delta')$ ($h$ states for the standard ruling). Such a family
$\Omega_N$ is found in \cite{\refMP}, Proposition 5.1B. A missing
weakly symplectic isotopy between $(C'_0,\omega'_0)$ and
$(C'_0,\omega'_0+N\pr^*\omega_{S^2})$, as well as that between $(
H_1,\omega'_1+N\pr^*\omega_{S^2})$ and $(
H_1,\omega'_1)$, can be given by variation of the symplectic
structure only: $\theta\mapsto\omega'_i+\theta\pr^*\omega_{S^2}$,
$i=1,2$.

The existence of $h_t$ with the properties enumerated above can
be proven in the following way.
Let number
the pseudo-lines going through $p_0$
by points in $S^2$, which
we identify with the projective
line of complex directions at $p_0$.
Then, to construct a desired $h_t$ we
fix a parametric representation of a pseudo-line
going through $p_0$ as follows.
We choose its intersection with a selected
pseudo-line not going through $p_0$ as the value at
$\infty$, $p_0$ as the value at $0$ and
fix the derivative of the parametrization
at $0$ (a unit tangent vector $\xi$
to $\C P^2$ at $p_0$). By means of such parametrizations,
$\phi_{\xi, t}:\C P^1\to\C P^2$,
we introduce the fibers of $h_t$ replacing $\phi_{\xi,t}$
in a disc $0\in D_1(\epsilon)\subset\C\subset\C P^1$
by $(1-\delta(r))\Xi+\delta(r)
\phi_{\xi,t}: D_1(\epsilon)\to
B(\delta)$, where $\Xi$
is the linear part of $\phi_{\xi,t}$ at $0$ and
$\delta:[0,\epsilon]\to [0,1]$ is a bump function which is taking value
$0$ near $0$ and $1$ near $\epsilon$.
The resulting maps $h_t$
have the enumerated above properties if $\epsilon$
is sufficiently small and $r\delta'(r)<1$
for any $r\in [0,\epsilon]$.
(One can check that a fibre in direction $\xi$ is symplectic by means
of straightforward calculations with $i dh\wedge\overline{dh}\wedge\Omega$
in affine complex coordinates
$z,w$ near $p_0$
with
$\Re \frac{\partial}{\partial z}=\xi$,
$\Re$ states for the real part,
using the complex analyticity
of the partical derivatives with
respect to local coordinates $s$ in $S^2$
of equation
$w=sz+\phi_2(s)z^2+\dots$
defining the pseudo-lines.
the calculations with $h$ one can use the above equation
to get an implicit equation for $h$:
$h=v-\delta(r)[\phi_2(h)z+\dots], w=vz$.)
\qed\enddemo

\proclaim{Corollary 4.3} A nodal symplectic, with respect to the
Fubini-Studi symplectic structure, surface is symplectically
$C^1$-smoothly isotopic to an algebraic curve if and only if its
braid factorization type $bmt$ with respect to a generic pencil is
a partial re-degeneration of some element from $S_{A_1}$:
$bmt=r(z_1)\cdot z_2$, where $r:S_{A_1}\to S_{A_0}$ is the
re-degeneration (see Example 1 in section 1.1), $A_0$ is the full
set of conjugates of the generator $a_1\in B_m$, $A_1$ is the full
set of conjugates of $a_1^2$, and $z_1,z_2\in S_{A_1}$.
\endproclaim

\demo{Proof} First, notice that $\alpha(bmt)=\Delta_m^2$, where
$m$ is the degree of the symplectic surface. So, according to
Theorem 1.2, $ z_1\cdot z_2=\widetilde\delta^2_m$. The latter
element is the braid monodromy factorization of $m$ lines in
general position. It remains to smooth the corresponding nodes of
this algebraic curve without deforming the other nodes, which is
possible, for example, by Bruzotti theorem \cite{\refBr}.
\qed\enddemo

\end

\end

\end

\end